\documentclass[12pt,leqno]{article}
\usepackage{amsfonts}
\usepackage{amsmath, amsthm, amsfonts, amssymb, color}
\usepackage{mathrsfs}
\usepackage{color}
\usepackage[notcite,notref]{showkeys}

\theoremstyle{definition}

\newcommand{\scr}[1]{\mathscr #1}
\definecolor{wco}{rgb}{0.5,0.2,0.3}

\numberwithin{equation}{section} \theoremstyle{remark}

\newcommand{\ua}{\uparrow}

\title{{\bf Bismut Formula for Lions Derivative of Distribution Dependent SDEs and Applications}\footnote{Supported in
 part by  NNSFC (11771326, 11831014, 11431014,11726627).} }
\author{
{\bf Panpan Ren$^{b)}$,    Feng-Yu Wang$^{a,b)}$  }\\
\footnotesize{$^{a)}$ Center for Applied Mathematics, Tianjin University, Tianjin 300072, China}\\
 \footnotesize{$^{b)}$ Department of Mathematics,
Swansea University, Singleton Park, SA2 8PP, United Kingdom}\\
\footnotesize{  673788@swansea.ac.uk, wangfy@tju.edu.cn, F.-Y.Wang@swansea.ac.uk}}
\begin{document}
\allowdisplaybreaks
\def\R{\mathbb R}  \def\ff{\frac} \def\ss{\sqrt} \def\B{\mathbf
B}
\def\N{\mathbb N} \def\kk{\kappa} \def\m{{\bf m}}
\def\ee{\varepsilon}\def\ddd{D^*}
\def\dd{\delta} \def\DD{\Delta} \def\vv{\varepsilon} \def\rr{\rho}
\def\<{\langle} \def\>{\rangle} \def\GG{\Gamma} \def\gg{\gamma}
  \def\nn{\nabla} \def\pp{\partial} \def\E{\mathbb E}
\def\d{\text{\rm{d}}} \def\bb{\beta} \def\aa{\alpha} \def\D{\scr D}
  \def\si{\sigma} \def\ess{\text{\rm{ess}}}
\def\beg{\begin} \def\beq{\begin{equation}}  \def\F{\scr F}
\def\Ric{\mathcal Ric} \def\Hess{\text{\rm{Hess}}}
\def\e{\text{\rm{e}}} \def\ua{\underline a} \def\OO{\Omega}  \def\oo{\omega}
 \def\tt{\tilde}
\def\cut{\text{\rm{cut}}} \def\P{\mathbb P} \def\ifn{I_n(f^{\bigotimes n})}
\def\C{\scr C}      \def\aaa{\mathbf{r}}     \def\r{r}
\def\gap{\text{\rm{gap}}} \def\prr{\pi_{{\bf m},\varrho}}  \def\r{\mathbf r}
\def\Z{\mathbb Z} \def\vrr{\varrho} \def\ll{\lambda}
\def\L{\scr L}\def\Tt{\tt} \def\TT{\tt}\def\II{\mathbb I}
\def\i{{\rm in}}\def\Sect{{\rm Sect}}  \def\H{\mathbb H}
\def\M{\scr M}\def\Q{\mathbb Q} \def\texto{\text{o}} \def\LL{\Lambda}
\def\Rank{{\rm Rank}} \def\B{\scr B} \def\i{{\rm i}} \def\HR{\hat{\R}^d}
\def\to{\rightarrow}\def\l{\ell}\def\iint{\int}
\def\EE{\scr E}\def\Cut{{\rm Cut}}\def\W{\mathbb W}
\def\A{\scr A} \def\Lip{{\rm Lip}}\def\S{\mathbb S}
\def\BB{\scr B}\def\Ent{{\rm Ent}} \def\i{{\rm i}}\def\itparallel{{\it\parallel}}
\def\g{{\mathbf g}}\def\Sect{{\mathcal Sec}}\def\T{\mathcal T}
\def\f{\mathbf f} \def\g{\mathbf g}\def\Zeta{\Xi}
\maketitle

\begin{abstract}By using Malliavin calculus, Bismut type formulas  are established for the Lions derivative   of $P_tf(\mu):=\E f(X_t^\mu)$, where $t>0,$ $ f $ is a bounded measurable function,  and $X_t^\mu$ solves a distribution dependent SDE with initial distribution $\mu$. As applications, explicit estimates are derived for the  Lions derivative and  the   total variational distance between distributions of   solutions  with different initial data. Both degenerate and non-degenerate situations are considered. Due to the lack of the semigroup property   and the invalidity of the formula $P_tf (\mu)= \int P_tf(x)\mu(\d x)$, essential difficulties are overcome in the study.
\end{abstract} \noindent
 \vskip 2cm

\section{Introduction}

The Bismut  formula  introduced in \cite{BB}, also called Bismut-Elworthy-Li formula due to \cite{EL},  is a powerful tool in characterising the regularity of distribution for SDEs and SPDEs. A plenty of results have been derived for this type formulas  and applications by using stochastic analysis  and coupling methods, see for instance \cite{Wbook} and references therein.

On the other hand, because of crucial applications in the study of  nonlinear PDEs and environment dependent financial systems, the distribution dependent SDEs (also called McKean-Vlasov or mean field SDEs) have received increasing attentions, see  \cite{DV1, DV2, FG, Gu, MV, SZ,V2} and references therein. Recently, this type SDEs have been applied in \cite{LP, HSS, LiJ, RW18} to characterize PDEs involving the Lions derivative ($L$-derivative for short) introduced  by P.-L. Lions in his lectures \cite{Card}. Moreover, Harnack inequality,   gradient estimates and exponential ergodicity have been investigated  in \cite{W18} and \cite{S18}.
In this paper, we aim to establish Bismut type $L$-derivative formula     for distribution dependent  SDEs with possibly degenerate noise.

To introduce our main results, we first recall the $L$-derivative.
Let $\scr P(\R^d)$ be the space of all probability measures on $\R^d$, and let
$$\scr P_2(\R^d)=\bigg\{\mu\in \scr P(\R^d):\ \mu(|\cdot|^2):=\int_{\R^d} |x|^2\mu(\d x)<\infty\bigg\}.$$ Then $\scr P_2(\R^d)$ is a Polish space under the Wasserstein distance
 $$\W_2(\mu,\nu):= \inf_{\pi\in \C(\mu,\nu)} \bigg(\int_{\R^d\times\R^d} |x-y|^2\pi(\d x,\d y)\bigg)^{\ff 1 2},\ \ \mu,\nu\in \scr P_2(\R^d),$$
 where $\C(\mu,\nu)$ is the set of couplings for $\mu$ and $\nu$; that is, $\pi\in\C(\mu,\nu)$ is a probability measure on $\R^d\times\R^d$ such that
 $\pi(\cdot\times\R^d)=\mu$ and $\pi(\R^d\times\cdot)=\nu$. We will use ${\bf 0}$ to denote vectors with components $0$, or the constant map taking value {\bf 0}.

\beg{defn}\label{defn1} Let  $f: \scr P_2(\R^d)\to \R$, and let $g: M\times \scr P_2(\R^d)\to \R$ for a differentiable manifold $M$.
\beg{enumerate} \item[(1)] $f$ is called  $L$-differentiable at $\mu\in \scr P_2(\R^d)$, if  the functional
$$L^2(\R^d\to \R^d,\mu)\ni \phi \mapsto f(\mu\circ({\rm Id} +\phi)^{-1})$$ is Fr\'echet differentiable at ${\bf0}\in L^2(\R^d\to \R^d,\mu)$;
that is, there exists $($hence, unique$)$ $\gg\in L^2(\R^d\to\R^d,\mu)$ such that
\beq\label{*D1}\lim_{\mu(|\phi|^2)\to 0} \ff{f(\mu\circ ({\rm Id}+\phi)^{-1})-f(\mu)- \mu(\<\gg,\phi\>)}{\ss{\mu(|\phi|^2)}}= 0.\end{equation}
In this case, we denote $D^L f(\mu)=\gg$ and call it the $L$-derivative of $f$ at $\mu$.
\item[(2)] If the $L$-derivative $D^L f(\mu)$ exists for all $\mu\in\scr P_2(\R^d)$, then $f$ is called $L$-differentiable. If, moreover, for every $\mu\in\scr P_2(\R^d)$
there exists a $\mu$-version $D^Lf(\mu)(\cdot)$ such that   $D^Lf(\mu)(x)$ is jointly continuous in  $(x,\mu)\in\R^d\times \scr P_2(\R^d)$, we denote $f\in C^{(1,0)}(\scr P_2(\R^d))$.
\item[(3)]  $g$   is called differentiable on $M\times \scr P_2(\R^d)$, if for any $(x,\mu)\in M\times \scr P_2(\R^d)$, $g(\cdot,\mu)$ is differentiable at $x$ and $g(x,\cdot)$ is $L$-differentiable at $\mu.$  If, moreover,
$\nn g(\cdot,\mu)(x)$   and  $D^L g(x,\cdot)(\mu)(y)$  are  joint continuous in $(x,y,\mu)\in M^2\times \scr P_2(\R^d)$, where $\nn$ is the gradient operator on $M$, we write $g\in C^{1,(1,0)} (M\times \scr P_2(\R^d))$. \end{enumerate}   \end{defn}
As indicated in \cite{RW18} that for any $n\ge 1$, $g\in C^1(\R^n)$ and $h_1,\cdots, h_n\in C_b^1(\R^d)$, the cylindrical function
$$\mu \mapsto g(\mu(h_1),\cdots, \mu(h_n))$$ is in $C^{(1,0)}(\scr P_2(\R^d))$ with
$$D^Lg(\mu)(x)= \sum_{i=1}^n \big(\pp_i  g(\mu(h_1),\cdots, \mu(h_n)) \big)\nn h_i(x),\ \ (x,\mu)\in \R^d\times \scr P_2(\R^d).$$

Obviously, if $f$ is $L$-differentiable at $\mu$, then
\beq\label{XMM1}   D_\phi^L f(\mu):= \lim_{\vv\downarrow 0} \ff{f(\mu\circ({\rm Id} +\vv\phi)^{-1})-f(\mu)}{\vv}
= \mu\big(\<D^L f(\mu), \phi\>\big), \ \phi\in L^2(\R^d\to\R^d,\mu). \end{equation}
We may call $D_\phi^L$ the directional $L$-derivative   along $\phi$, which was introduced in \cite{AKR, OTT}.

When $D_\phi^Lf(\mu)$  is a bounded linear functional of $\phi\in L^2(\R^d\to\R^d,\mu)$, there exists a unique $\xi\in L^2(\R^d\to\R^d,\mu)$ such that
$D_\phi^Lf(\mu)= \mu(\<\xi, \phi\>)$ holds for all $\phi\in L^2(\R^d\to\R^d,\mu)$. In this case, $\phi\mapsto f(\mu\circ({\rm Id}+\phi)^{-1} )$ is G\^ateaux differentiable at ${\bf 0}$, and we say that $f$ is weakly $L$-differentiable at $\mu$, since the G\^ateaux differentiability is weaker than the Fr\'echet one.

By \eqref{XMM1}, for an $L$-differentiable function $f$ on $\scr P_2(\R^d)$, we have
\beq\label{XMM} \|D^Lf(\mu)\|:= \|D^Lf(\mu)(\cdot)\|_{L^2(\mu)}=\sup_{\mu(|\phi|^2)\le 1} |D_\phi^L f(\mu)|.\end{equation}
For a vector-valued  function $f=(f_i)$, or  a matrix-valued function $f=(f_{ij})$  with $L$-differentiable components, we write
$$D_\phi^L f(\mu)= (D_\phi^L f_i(\mu)),\ \text{or}\  D_\phi^Lf(\mu)=(D_\phi^L f_{ij}(\mu)),\ \ \ \mu\in \scr P_2(\R^d).$$

 Let $W_t$ be a $d$-dimensional Brownian motion on the natural filtered probability space $(\OO^0,\F^0,\{\F_t^0\}_{t\ge 0},\P)$. To ensure that for any $\mu\in \scr P_2(\R^d)$ there exists a random variable $X$ on $\R^d$ with distribution $\mu$,   let $\mu^0$  be a probability   measure on $\R^d$ which is equivalent to the Lebesgue measure, and enlarge the probability space as
 $$(\OO,\F, \{\F_t\}_{t\ge 0},\P):= (\OO^0\times\R^d, \F^0\times \B(\R^d), \{\F_t^0\times \B(\R^d)\}_{t\ge 0},\ \P^0\times\mu^0).$$ Then
 $$W_t(\oo):= W_t(\oo^0),\ \ t\ge 0, \oo:=(\oo^0,x) \in \OO$$ is a $d$-dimensional Brownian motion on $(\OO,\F, \{\F_t\}_{t\ge 0},\P).$
Let $\L_\xi$ denote the distribution of a random variable on the probability space $(\OO,\F,\P)$. In case different probability spaces are concerned, we write $\L_{\xi|\P}$ instead of $\L_\xi$ to emphasize the reference probability measure $\P$.

Consider the following distribution dependent SDE on $\R^d$:
\beq\label{E1} \d X_t=  b_t(X_t,\L_{X_t})\d t + \si_t(X_t,\L_{X_t}) \d W_t,\ \ X_0\in L^2(\OO\to \R^d,\F_0,\P),\end{equation}
where
$$\si: [0,\infty)\times \R^d\times \scr P_2(\R^d)\to \R^{d\otimes d},\ \ b: [0,\infty)\times \R^d\times \scr P_2(\R^d)\to \R^d$$
are continuous such that for some increasing function $K:\, [0,\infty)\to [0,\infty)$ there holds
\beq\label{SPP1}
\beg{split} &|  b_t(x,\mu)- b_t(y,\nu)|+\|\si_t(x,\mu)-\si_t(y,\nu)\|\\
&\le K(t)\big(|x-y|+\W_2(\mu,\nu)\big),\ \ t\ge 0, x,y\in \R^d, \mu,\nu\in \scr P_2(\R^d)
\end{split}\end{equation}
 and
 \beq\label{SPPn}
\|\si_t({\bf 0},\dd_{{\bf 0}})\|+ |b_t({\bf 0},\dd_{{\bf 0}})|\le K(t),\ \ t\ge 0,
\end{equation}
where and in what follows, for $x\in\R^d$ we denote  by $\dd_{x}$   the Dirac measure at $x$, and $\|\cdot\|$ is the operator norm. For any $t\ge 0$, let $L^2(\OO\to\R^d,\F_t,\P)$ be the class of $\F_t$-measurable square integrable random variables on $\R^d$. By \eqref{SPP1} and \eqref{SPPn},
for any $s\ge 0$ and $X_s\in L^2(\OO\to\R^d,\F_s,\P)$, \eqref{E1} has a unique solution $(X_{s,t})_{t\ge s}$ with $X_{s,s}=X_s$ and
\beq\label{MM}
 \E \Big[\sup_{t\in [s,T]} |X_{s,t}|^2\Big]<\infty,\ \ T\ge s,
\end{equation} see, for instance
  \cite{W18},  where   gradient estimates and Harnack inequalities are also derived for the associated nonlinear semigroup. See also \cite{18HW,MV}   for weaker conditions ensuring  the existence and uniqueness of solutions to \eqref{E1}.
For any $\mu\in \scr P_2(\R^d)$ and $s\ge 0$,       let $(X_{s,t}^\mu)_{t\ge s}$ be the solution to \eqref{E1} with $\L_{X_{s,s}}= \mu.$   Denote
\beq\label{SD1} P_{s,t}^*\mu  =\L_{X_{s,t}^\mu},\ \    t\ge s, \mu\in \scr P_2(\R^d).\end{equation}  Let
\beq\label{SD2} (P_{s,t}f)(\mu)=(P_{s,t}^*\mu)(f):= \int_{\R^d} f\d (P_{s,t}^*\mu)= \E f(X_{s,t}^\mu),\ \ t\ge s, f\in \B_b(\R^d), \mu\in\scr P_2(\R^d).\end{equation} Then  for any $0\le s\le t$, $P_{s,t}$ is a linear operator from   $\B_b(\R^d)$ to $\B_b(\scr P_2(\R^d))$.

In this paper, we aim to establish the Bismut type formula for the  $L$-derivative of $P_{s,t}f$ for $t>s$. By considering the SDE for
$\tt X_t:= X_{t+s}, t\ge 0$,   without loss of generality we may and do assume $s=0$. So, for simplicity, below we only establish the derivative formula for $P_tf:= P_{0,t}f, t>0$.
  More precisely, for any $T>0$, $\mu\in\scr P_2(\R^d)$ and $\phi\in L^2(\R^d\to \R^d,\mu)$,  we aim
to construct an integrable  random variable $M_T^{\mu,\phi}$ such that
\beq\label{BSM} D^L_\phi (P_T f)(\mu)= \E \big[f(X_T^\mu)M_T^{\mu,\phi}\big],\ \ f\in \B_b(\R^d),\end{equation}  which in turn  implies the $L$-differentiability of $P_Tf$.
 Note that   the derivative formula for $(P_Tf)(x):= (P_Tf)(\dd_x)$ along a vector $v\in \R^d$ is derived in \cite{BAN}, which is the special case of \eqref{BSM} with $\mu=\dd_x$ and $\phi\equiv v$.
Moreover,  formulas of  the $L$-derivative and integration by parts  have been  presented in \cite{CM} for the following de-coupled SDE:
$$\d X_t^{x,\mu}=   b(t, X_t^{x,\mu}, P_t^*\mu)\d t + \si(t,X_t^{x,\mu}, P_t^*\mu) \d W_t, \ \ X_0^{x,\mu}=x,$$
which is different from the original SDE \eqref{E1} but has important applications in solving PDEs with Lions' derivatives, see \cite{LP, LiJ, RW18} and references within.

When the SDE \eqref{E1} is distribution independent, i.e. $b_t(x,\mu)= b_t(x)$ and  $\si_t(x,\mu)=\si_t(x)$ do not depend on $\mu$, the Bismut type formula
\beq\label{BSM0} \nn P_T f(x)= \E\big[f(X_T^x) M_T^x\big],\ \   x\in\R^d, f\in \B_b(\R^d)\end{equation}  has been well studied in the literature, where $M_T^x$ is an integrable random variable on $\R^d$, which is measurable in $x\in \R^d$ when it varies, see for instance  \cite{TH, GW12,   W14, W16, WZ13} and references within.   Since the coefficients are distribution independent, we have
\beq\label{NLL}(P_Tf)(\mu) =\int_{\R^d} (P_T f)(x) \,\mu(\d x),\ \ \end{equation}
so that $P_Tf$ is $L$-differentiable with $D^L (P_Tf)(\mu)=\nn P_T f.$ Hence, by \eqref{BSM0} and \eqref{NLL} we obtain
\beg{align*} D^L_\phi (P_T f)(\mu) &= \mu(\<D^LP_T f, \phi\>)  = \int_{\R^d} \E\big[f(X_T^x)\<M_T^x, \phi(x)\>\big]\mu(\d x) \\
 &= \E\big[f(X_T^\mu)\<M_T^{X_0^\mu},\phi(X_0^\mu)\>\big].\end{align*} Therefore,     \eqref{BSM} holds for   $M_T^{\mu,\phi}=\<M_T^{X_0^\mu},\phi(X_0^\mu)\>$.

However, when the SDE is distribution dependent,
as explained in \cite{W18} that    in general
\eqref{NLL} does  not hold, so  it is non-trivial to establish the Bismut type formula \eqref{BSM}.

\

The remainder of the paper is organized as follows. In section 2, we state our main results on Bismut formulas  of  $D_\phi^LP_Tf$  and applications,  for both non-degenerate and  degenerate distribution dependent SDEs.
To establish the Bismut formula   using Malliavin calculus, we make necessary preparations in Section 3 concerning partial derivatives in the initial value, and   Malliavin derivative   for solutions of \eqref{E1}.  Finally,  complete proofs of the main results are addressed in Section 4.

\section{Main results}

 Let $|\cdot|$ denote the Euclidean norm in $\R^d$, and $\|\cdot\|$ denote the operator norm for matrices or more generally linear operators. We   make the following assumption.

\beg{enumerate}\item[{\bf (H)}] For any $t\ge 0$, $b_t,\si_t\in C^{1,(1,0)}(\R^d\times \scr P_2(\R^d))$. Moreover,  there exists a continuous function $K:[0,\infty)\to [0,\infty),$ such that   \eqref{SPPn} holds and for any $\mu\in\scr \P_2(\R^d)$, 
\beg{align*} & \|\nn b_t(\cdot,\mu)(x)\|,  \|D^L b_t(x,\cdot)(\mu)\|_\infty,\  \ff 1 2\|\nn \si_t(\cdot,\mu)(x)\|^2,  \ff 1 2 \|D^L \si_t(x,\cdot)(\mu)\|^2 _\infty \Big\}\le K_t,\\
&  |D^L b_t(x,\cdot)(\mu)(y)|+  |D^L \si_t(x,\cdot)(\mu)(y)|^2\le K_t(1+|x|^2 +|y|^2),\  \  t\ge 0, x\in\R^d,  \end{align*} where   $\|D^Lf(\mu)\|_\infty:= \|D^Lf(\mu)(\cdot)\|_{L^\infty(\mu)}$ for an $L$-differentiable function $f$ at $\mu$. Moreover, 
there exists a 
  \end{enumerate}
Obviously, {\bf (H)} implies \eqref{SPP1} and \eqref{SPPn}, so that the SDE \eqref{E1} has a unique solution for any  initial value $X_0\in L^2(\OO\to\R^d,\F_0,\P)$.

In the following two subsections, we state our main results for non-degenerate and degenerate cases respectively.

\subsection{The non-degenerate case}
Due to technical reasons,  the following result Theorem \ref{T3.1} only works  for distribution independent $\si_t$. But  some other results (for instance
Proposition \ref{P2.1}) apply  to the   general setting. So, in addition to {\bf (H)}   we also assume
 \beq\label{LLP} \si_t(x,\mu)=\si_t(x)\ \text{ with\ }\|\si_t(x)^{-1}\|\le \ll_t \ \text{for\ some\ }   \ll\in C([0,\infty)\to(0,\infty)).\end{equation}  Let $\mu\in \scr P_2(\R^d)$,  and let $X_t$ solve \eqref{E1} for $X_0 \in L^2(\OO\to\R^d,\F_0,\P)$ with $\L_{X_0}=\mu$.
Given $\phi\in L^2(\R^d\to\R^d,\mu)$, consider the following  SDE for $v_t^\phi$ on $\R^d$:
\beq\label{DR0} \beg{split} \d v_t^\phi & =\Big\{\nn_{v_t^\phi} b_t(\cdot, \L_{X_t})(X_t)+ \big(\E \< D^L b_t(y,\cdot)(\L_{X_t})(X_t), v_t^\phi\>\big)\big|_{y=X_t}\Big\}\d t\\
&+ \Big\{\nn_{v_t^\phi} \si_t(X_t)\Big\}\d W_t,\ \
v_0^\phi=\phi(X_0).\end{split}\end{equation} By {\bf (H)}, this linear SDE is   well-posed  with $\sup_{t\in [0,T]}\E |v_t^\phi|^2\le C \mu(|\phi|^2)$  for some constant $C=C(T)>0$, see \eqref{PPG} below.
Denote  $g_s'= \ff{\d}{\d s} g_s$ for a differentiable function $g$ of $s\in \R$.

\beg{thm}\label{T3.1} Assume {\bf (H)} and $\eqref{LLP}$.
Then for any   $f\in \B_b(\R^d),\mu\in\scr P_2(\R^d)$ and $T>0$, $P_Tf$ is  $L$-differentiable at $\mu$ such that  for any $g\in C^1([0,T])$ with $g_0=0$ and $g_T=1$,
\beq\label{BSM1} \beg{split} D_\phi^L(P_Tf)(\mu) &= \E\bigg[f(X_T)   \int_0^T \big\< \zeta_t^\phi,\ \d W_t\big\>\bigg], \  \phi\in L^2(\R^d\to\R^d,\mu),\end{split}\end{equation} where   $X_t$ solves $\eqref{E1}$ for $\L_{X_0}=\mu$, and
$$  \zeta_t^\phi:= \si_t(X_t)^{-1}\Big\{ g_t' v_t^\phi+ \big(\E \< D^L b_t(y,\cdot)(\L_{X_t})(X_t), g_t v_t^\phi \>\big)\big|_{y=X_t}\Big\},\ \ t\in [0,T].$$ Moreover, the limit
\beq\label{*D*} D_\phi^LP_T^*\mu:= \lim_{\vv\downarrow 0} \ff{P_T^*\mu\circ({\rm Id}+\vv\phi)^{-1}- P_T^*\mu}\vv= \psi P_T^*\mu\end{equation}  exists in the total variational norm, where
  $\psi$ is  the unique element in $L^2(\R^d\to\R, P_T^*\mu) $ such that
$\psi(X_T)=  \E\big(\int_0^T \big\<\zeta_t^\phi,\ \d W_t\big\>\big|X_T\big),$ and $(\psi P_T^*\mu)(A):= \int_A \psi\d P_T^*\mu,\ A\in\B(\R^d)$.
  \end{thm}

\paragraph{Remark 2.1.} When $f\in C_b^1(\R^d)$, \eqref{BSM1} can be proved  as in the distribution independent case by constructing  a proper random variable  $h$ on the Cameron-Martin space such that   $D_h X_T=\nn_{\phi}X_T$. However,  for the $L$-differentiability of $P_Tf$, one has to construct $\gg\in L^2(\R^d\to\R^d,\mu)$ such that \eqref{*D1} holds for $P_Tf$ replacing $f$, which  is non-trivial.

Moreover, comparing with the classical case where \eqref{BSM1} for $f\in C_b^1(\R^d)$ can be easily extended to $f\in \B_b(\R^d)$, there is essential difficulty to do this in the   distribution dependent setting. More precisely, when $b_t$ and $\si_t$ do not depend on the distribution,   we have the semigroup property
  $P_Tf(\mu)= P_t(P_{t,T} f)(\mu)$ for $t\in (0,T)$, where $P_{t,T}f(x):= P_{t,T}f(\dd_x) $ for the Dirac measure $\dd_x$ at point $x$. In many cases, we have  $P_{t,T}f \in C_b^1(\R^d)$ for $f\in \B_b(\R^d)$.  Then for any $f\in \B_b(\R^d)$,  one may apply the derivative formula \eqref{BSM1} with $(P_{t},P_{t,T}f)$ replacing $(P_T,f)$ to derive a derivative formula for $P_Tf$. However, in the distribution dependent case, due to the lack of \eqref{NLL}  we no longer have $P_Tf(\mu)= P_t(P_{t,T}f)(\mu)$, so that this argument becomes invalid. To overcome this difficulty we will make a new approximation argument, see    step (a) in the proof of Theorem \ref{T3.1} for details.

\

As applications of Theorem \ref{T3.1}, the following result  consists of   estimates   on the $L$-derivative and   the total variational distance between distributions of  solutions with different initial data. 

\beg{cor} \label{T3.2} Assume {\bf (H)} and $\eqref{LLP}$ for some increasing functions $K$ and continuous function $\ll$.
   \beg{enumerate}\item[$(1)$] For any $f\in \B_b(\R^d)$ and $T>0$,
\beq\label{GRD} \beg{split}&\| D^L (P_Tf)(\mu)\|^2:= \sup_{\mu(|\phi|^2)\le 1 } |D^L_\phi (P_Tf)(\mu)|^2 \\
&\le \big\{(P_Tf^2)(\mu)-(P_t f(\mu))^2\big\} \int_0^T \Big(\ff 1 T+ K_t\Big)^2  \ll_t^{2} \e^{8K_tt}\d t.\end{split}\end{equation}
\item[$(2)$]  For any $T>0$, \beq\label{GRD3}\beg{split} & |P_Tf(\mu)- P_Tf(\nu)|^2\\
&\le   \|f\|_\infty^2 \W_2(\mu,\nu)^2\int_0^T \Big(\ff 1 T+ K_t\Big)^2  \ll_t^{2} \e^{8K(t)t}\d t,\ \   \mu,\nu\in \scr P_2(\R^d), f\in \B_b(\R^d).\end{split}\end{equation} Consequently, for any $T>0$ and $\mu,\nu\in \scr P_2(\R^d)$,
\beq\label{GRD2} \beg{split}&\|P_T^*\mu-P_T^*\nu\|_{var}^2:= \sup_{A\in \B(\R^d)} |(P_T^*\mu)(A)-(P_T^*\nu)(A)|^2\\
 &\le  \W_2(\mu,\nu)^2\int_0^T \Big(\ff 1 T+ K_t\Big)^2  \ll_t^{2} \e^{8K(t)t}\d t.\end{split}\end{equation}
  \end{enumerate}
 \end{cor}

\subsection{Stochastic Hamiltonian systems}

Consider the following distribution dependent stochastic Hamiltonian system for $X_t=(X_t^{(1)}, X_t^{(2)})$ on $\R^{m+d}=\R^m\times \R^d$:
\beq\label{E5} \beg{cases} \d X_t^{(1)}= b^{(1)}_t(X_t)\d t,\\
\d X_t^{(2)} = b_t^{(2)}(X_t, \L_{X_t})\d t +\si_t \d W_t,\end{cases}\end{equation} where
  $(W_t)_{t\ge 0}$ is a $d$-dimensional Brownian motion as before, and  for each $t\ge 0$, $\si_t$ is an invertible $d\times d$-matrix,
$$b_t= (b_t^{(1)}, b_t^{(2)}): \R^{m+d}\times \scr P_2(\R^{m+d}) \to  \R^{m+d}$$ is measurable with $b^{(1)}_t(x,\mu)= b_t^{(1)}(x)$ independent of the distribution $\mu$.
Let $\nn=(\nn^{(1)}, \nn^{(2)})$ be the gradient operator on $\R^{m+d}=\R^m\times\R^d$, where $\nn^{(i)}$ is the gradient in the $i$-th component, $i=1,2$.
Let $\nn^2=\nn\nn$ denote the Hessian operator on $\R^{m+d}$. We assume

\beg{enumerate}\item[{\bf (H1)}] For every $t\ge 0$, $b_t^{(1)} \in C^2_b(\R^{m+d}\to\R^m),$ $b_t^{(2)}\in C^{1,(1,0)}(\R^{m+d}\times \scr P_2(\R^{m+d})\to\R^d)$, and there exists an increasing function $K: [0,\infty)\to [0,\infty)$ such that \eqref{SPPn} and
$$\|\nn b_t(\cdot,\mu)(x)\|+ \|D^L b_t^{(2)}(x,\cdot)(\mu)\|  +\|\nn^2 b_t^{(1)}(\cdot,\mu)(x)\|  \le K(t)$$ hold for all $ t\ge 0, (x,\mu)\in\R^d\times \scr P_2(\R^d)$.
  \end{enumerate}

Obviously,  this assumption implies {\bf (H)} for   the SDE \eqref{E5}.
We aim to establish the derivative formula of type \eqref{BSM} with $P_t$ and $P_t^*$  being  defined by \eqref{SD1} and \eqref{SD2} for the SDE \eqref{E5}. To  follow the line of \cite{WZ13} where the distribution independent model was investigated,
  we need the following assumption {\bf (H2)}.

   For any $s\ge 0,$ let $\{K_{t,s}\}_{t\ge s}$ solve the following linear random ODE on $\R^{m\otimes m}$:
\begin{align}
\ff{\d}{\d t}K_{t,s}=  (\nn^{(1)}b^{(1)})(X_t) K_{t,s},\ \ \ t\ge s,  K_{s,s}=I_{m\times m},\label{Eq1}
\end{align} where $I_{m\times m}$ is the $m\times m$-order identity matrix.

\beg{enumerate} \item[{\bf (H2)}] There exists  $B\in \B_b([0,T]\to \R^{m\otimes d}) $ such that
 \beq\label{B} \<(\nn^{(2)} b^{(1)}_t - B_t)B_t^* a,a\>\ge -\vv |B_t^*a|^2,\ \ \forall a\in \R^m\end{equation}   holds for some constant $\vv\in [0,1).$
Moreover,  there exists an increasing function $\theta\in C([0,T])$ with $\theta_t>0$ for $t\in (0,T]$
such that
\beq\label{B2} \int_0^t s(T-s) K_{T,s}B_sB_s^*K_{T,s}^*\d s\ge \theta_t  I_{m\times m},\ \ t\in (0,T].\end{equation}
\end{enumerate}

\paragraph{Example 2.1.}  Let  $$b_t^{(1)}(x)= A x^{(1)} + Bx^{(2)},\ \  x=(x^{(1)}, x^{(2)}) \in \R^{m+d}$$ for some  $m\times m$-matrix $A$ and $m\times d$-matrix $B$. If
 the Kalman's rank condition
$$ \text{Rank}[B, AB, \cdots, A^{k}B]=m $$ holds for some $k\ge 1$, then  {\bf (H2)} is satisfied  with $\theta_t=c_T t$ for some constant $c_T>0$, see the proof of \cite[Theorem 4.2]{WZ13}.
In general, {\bf (H2)} remains true under small perturbations of this $b_t^{(1)}.$

\

According to the proof of \cite[Theorem 1.1]{WZ13}, {\bf (H2)} implies that  the matrices
$$Q_t:=\int_0^t s(T-s) K_{T,s}\nn^{(2)}b_s^{(1)}(X_s) B_s^* K_{T,s}^* \d s,\ \   t\in (0,T]$$ are invertible with
\beq\label{Q}\|Q_t^{-1}\|\le\ff 1 { (1-\vv)\theta_t },\ \ t\in (0,T].\end{equation} For $(X_t)_{t\in [0,T]}$ solving \eqref{E5} with $\L_{X_0}=\mu$ and $\phi=(\phi^{(1)},\phi^{(2)})\in L^2(\R^{m+d}\to\R^{m+d},\mu)$, let
\beg{equation}\label{aa}\beg{split} \aa_t^{(2)}  = &\ff{T-t} T \phi^{(2)} (X_0) -\ff{t(T-t) B_t^*K_{T,t}^*}{\int_0^T \theta_s^2 \d s} \int_t^T  \theta_s^2 Q_s^{-1} K_{T,0}\phi^{(1)}(X_0)\d s\\
&-t(T-t) B_t^* K_{T,t}^*Q_T^{-1}\int_0^T\ff{T-s}TK_{T,s}\nn^{(2)}_{\phi^{(2)} (X_0)} b^{(1)}_s(X_s)\d s,\ \ t\in [0,T],
\end{split}\end{equation} and
\beg{equation}\label{B0} \aa^{(1)}_t= K_{t,0}\phi^{(1)} (X_0)+\int_0^tK_{t,s}\nn^{(2)}_{\aa_s^{(2)}} b_s^{(1)}(X_s(x))\,\d s, \ \ t\in [0,T]. \end{equation}  Moreover, let $(h_t^\aa,w_t^\aa)_{t\in [0,T]}$ be the unique solution to   the random ODEs
\beq\label{B00}\beg{split}& \ff{\d h_{t}^\aa}{\d t} =  \si_t^{-1} \Big\{\nn_{\aa_t} b_t^{(2)}(X_t,\L_{X_t}) -(\aa_t^{(2)})'\\
&\qquad\qquad\qquad + \big(\E\<D^L b_t^{(2)}(y,\cdot)(\L_{X_t})(X_t), \aa_t+w_{t}^\aa\>\big)\big|_{y=X_t}\Big\},\\
&\ff{ \d w_{t}^\aa}{\d t} = \nn_{w_{t}^\aa} b_t (\cdot,\L_{X_t})(X_t) + ({\bf 0}, \si_t (h_{t}^\aa)'),\ \ \ \ h_0^\aa= w_0^\aa=0. \end{split}\end{equation}
 Let $ (D^*, \D(D^*))$ be the   Malliavin divergence operator associated with the Brownian motion $(W_t)_{t\in [0,T]}$, see Subsection 3.2 below for details. Then the main result in this part is the following.

\beg{thm}\label{T4.2} Assume {\bf (H1)} and {\bf (H2)}. Then $h^\aa\in \D(D^*)$ with    $\E|D^*(h^{\aa})|^p<\infty$ for all $p\in [1,\infty)$. Moreover, for any $f\in \B_b(\R^{m+d})$ and $T>0$, $P_Tf$ is $L$-differentiable at $\mu$ such that
\beq\label{BSMN} D_\phi^L(P_Tf)(\mu)= \E\big[f(X_T)\, D^*(h^\aa)\big].\end{equation}
Consequently:
\beg{enumerate} \item[$(1)$] $\eqref{*D*}$ holds for the unique $\psi\in L^2(\R^{m+d}\to\R, P_T^*\mu)$ such that $\psi(X_T)= \E(D^*(h^\aa)|X_T).$
\item[$(2)$] There exists  a  constant $c\ge 0$ such that for any $T>0$,
\beq\label{LD1}   \|D^L  (P_T f)(\mu)\|\le c\ss{ P_T|f|^2(\mu) -(P_Tf)^2(\mu)}  \ff{\ss{T}(T^2 +\theta_T )}{\int_0^{T} \theta_s^2\d s},\ \ f\in \B_b(\R^{m+d}),\end{equation}
\beq\label{LD2} \|P_T^*\mu- P_T^*\nu\|_{var} \le c \W_2(\mu,\nu) \ff{\ss{T}(T^2 +\theta_T) }{\int_0^{T} \theta_s^2\d s},\ \  \mu,\nu\in \scr P_2(\R^d).\end{equation}\end{enumerate}
 \end{thm}

\section{Preparations }
We first introduce a formula of the $L$-derivative re-organized from \cite[Theorem 6.5]{Card} and \cite[Proposition A.2]{HSS}, then
investigate  the partial derivatives  of $X_t$ in the initial value, and the Malliavin derivatives of $X_t$ with respect to the Brownian motion $W_t$.

\subsection{A formula of $L$-derivative}
The following result is essentially due to \cite[Theorem 6.5]{Card} for $f\in C^{(1,0)}(\scr P_2(\R^d))$, and   \cite[Proposition A.2]{HSS} for bounded $X$ and $Y$. We include a complete proof for readers' convenience.

\beg{prp}\label{P01} Let $(\OO,\F,\P)$ be an atomless probability space, and let $X,Y \in L^2(\OO\to\R^d,\P)$ with $\L_X=\mu$. If either $X$ and $Y$ are bounded and $f$ is $L$-differentiable at $\mu$, or $f\in C^{(1,0)}(\scr P_2(\R^d))$, then
\beq\label{XMM2}\lim_{\vv\downarrow 0} \ff{f(\L_{X+\vv Y})- f(\mu)}\vv = \E\<D^L f(\mu)(X), Y\>.\end{equation}
 Consequently,
 \beq\label{GMM} \Big|\lim_{\vv\downarrow 0} \ff{f(\L_{X+\vv Y})- f(\mu)}\vv \Big| =\big|\E\<D^L f(\mu)(X), Y\>\big|\le  \|D^Lf(\mu)\| \ss{\E |Y|^2}.\end{equation}
\end{prp}
\beg{proof}  It is easy  to see that   \eqref{GMM} follows from   \eqref{XMM} and \eqref{XMM2}. Indeed, letting $\phi\in L^2(\R^d\to\R^d,\mu)$ such that $\phi(X)= \E(Y|X)$, we have
\beg{align*} &\big|\E\<D^L f(\mu)(X), Y\>\big|=\big|\E\<D^L f(\mu)(X), \phi(X)\>\big| = \big|\mu(\<D^Lf(\mu), \phi\>)\big| \\
 &\le  \|D^Lf(\mu)\| \cdot \|\phi\|_{L^2(\mu)} = \|D^Lf(\mu)\| \big(\E|\E(Y|X)|^2\big)^{\ff 1 2} \le   \|D^Lf(\mu)\|\ss{\E |Y|^2}.\end{align*}
Below we prove \eqref{XMM2} for the stated two situations respectively.

(1) Assume that $X$ and $Y$ are bounded. For any $\R^d$-valued random variable $\xi$, let
$F(\xi)= f(\L_\xi).$ Next, let $(\bar\OO,\bar\F,\bar \P)$ be an atomless Polish probability space, and let $\bar X\in L^2(\bar\OO\to\R^d,\bar \P)$ with $\L_{\bar X|\bar \P }=\mu$, where $\L_{\cdot |\bar \P } $ denotes the distribution of a random variable under $\bar \P$. According to \cite[Proposition A.2(iii)]{HSS}, if
$$\bar F(\bar Y):= f(\L_{\bar Y|\bar \P}),\ \ \bar Y\in L^2(\bar\OO\to\R^d,\bar\P)$$
is Fr\'echet differentiable at $\bar X$ with derivative $D \bar F(\bar X)= D^Lf(\mu)(\bar X),$ then
 \beq\label{FY0} \lim_{\vv\downarrow 0} \ff{f(\L_{X+\vv Y})-f (\L_X)-\vv \E\<D^L  f(\mu)(X), Y\>}{\vv}=0.\end{equation}
Equivalently,  \eqref{XMM2} holds.
Below we construct the desired $\bar X$ and $(\bar\OO,\bar\F,\bar \P)$ such that $D \bar F(\bar X)= D^Lf(\mu)(\bar X).$

A natural choice of $(\bar\OO,\bar\F,\bar \P)$  is $(\R^d,\B(\R^d),\mu)$, but to ensure the  atomless property, we take $(\bar\OO,\bar\F,\bar \P)=(\R^d\times \R, \B(\R^d\times \R), \mu\times \ll),$ where $\ll$ is the standard Gaussian measure on $\R$.
Then $(\bar\OO,\bar\F,\bar \P)$ is an atomless Polish probability space. Let
$$\bar X(\bar\oo)= x,\ \ \bar\oo=(x,r)\in \R^d\times \R.$$
We have  $\L_{\bar X}=\mu$. 
Moreover, let
$$\tt f(\tt\mu)= f(\tt\mu(\cdot\times \R)),\ \ \tt\mu\in \scr P_2(\R^d\times \R).$$ It is easy to see that the $L$-differentiability of $f$ at $\mu$ implies that of $\tt f$ at $\mu\times\dd_0$ with
\beq\label{FY2} D^L\tt f (\mu\times \dd_0)(x,r)= (D^Lf(\mu)(x), 0),\ \ (x,r)\in \R^d\times\R.\end{equation}
Finally, on the probability space $(\OO,\F,\P)$ we have
\beq\label{FY1} F(Y):= f(\L_Y)= \tt f(\L_{\tt Y}),\ \ \tt Y:=(Y, 0)\in L^2(\OO\to \R^d\times\R,\F,\P).\end{equation}
Letting $\tt X=(X,0)\in L^2(\OO\to\T^d\times\R, \F,\P)$,    by \cite[Proposition A.2(iii)]{HSS}, the formula \eqref{FY0} holds for $(\tt X,\tt Y,\tt f, \mu\times\dd_0)$ replacing  $(X,Y, f, \mu)$, i.e.
\beg{align*}  \lim_{\vv\downarrow 0}   \ff{\tt f(\L_{\tt X+\vv\tt Y})-\tt f (\L_{\tt X})-\E\<D^L \tt f(\mu\times\dd_0),\vv \tt Y\>}{\vv} =0.\end{align*} Combining this with \eqref{FY2} and \eqref{FY1}, we prove \eqref{FY0}. Therefore,   \eqref{XMM2} holds.

(2) Let $f\in C^{(1,0)}(\scr P_2(\R^d))$ and  let $\mu\in \scr P_2(\R^d)$ and $X\in L^2(\OO\to\R^d,\P)$ with $\L_{X}=\mu$. For any $n\ge 1$, let
 $$x_n= \ff{x}{\ss{1+n^{-1}|x|^2}},\ \ x\in \R^d.$$ By \eqref{XMM2} for bounded $X$ and $Y$,  for any $n\ge 1$ we have
 \beq\label{FY3} \beg{split} & f(\L_{X_n+\vv Y_n})- f(\L_{X_n})= \int_0^\vv \ff{\d}{\d s} f(\L_{X_n+sY_n})\,\d s\\
 &= \int_0^\vv \E \<D^L f(\L_{X_n+sY_n})(X_n+sY_n), Y_n\>\,\d s.\end{split} \end{equation} Since $f\in C^{(1,0)}(\scr P_2(\R^d)),$  it follows that
$$  \sup_{n\ge 1,s\in [0,\vv]} \|D^L f(\L_{X_n + s Y_n})\|<\infty,\ \ \lim_{n\to\infty} \{ f(\L_{X_n+\vv Y_n})- f(\L_{X_n})\}=  f(\L_{X+\vv Y})- f(\L_{X}),$$ and for any $s\in [0,\vv]$,
$$ \lim_{n\to\infty} \E\big(|X-X_n|^2+|Y-Y_n|^2+ |D^L f(\L_{X_n+sY_n})(X_n+sY_n)- D^L f(\L_{X+sY})(X+sY)|^2 \big)=0.$$
 Then letting $n\to\infty$ in \eqref{FY3} we arrive at
\beq\label{QP1}  f(\L_{X+\vv Y})- f(\L_{X}) = \int_0^\vv \E \<D^L f(\L_{X+sY})(X+sY), Y\>\,\d s,\ \ \vv>0.\end{equation}
This implies \eqref{XMM2}. More precisely, it is easy to see that  $\{\L_{X+sY}\}$ is compact in $\scr P_2(\R^d)$. So,    $f\in C^{(1,0)}(\scr P_2(\R^d))$ implies
\beq\label{QP2} A:= \sup_{s\in [0,1]} \ss{\E|D^L f(\L_{X+sY})(X+sY)|^2} = \sup_{s\in [0,1]} \|D^Lf(\L_{X+sY})\|_{L^2(\L_{X+sY})}<\infty.\end{equation}
Combining this with the continuity property of $D^Lf$ on $\R^d\times\scr P_2(\R^d)$, we conclude that
$$\lim_{\vv\downarrow 0} D^Lf(\L_{X+sY})(X+sY) =D^L f(\L_X)(X)\  \text{weakly\ in\ } L^2(\OO\to\R^d,\P).$$ In particular,
\beq\label{QP3} \lim_{\vv\downarrow 0} \E\<D^L f(\L_{X+sY})(X+sY),Y\> =\E\<D^L f(\L_X) (X), Y\>.\end{equation} Moreover, \eqref{QP2} implies
$$\sup_{s\in [0,1]} \E\big|\<D^L f(\L_{X+sY})(X+sY),Y\>\big| \le A\ss{\E|Y|^2} <\infty.$$ Due to this,  \eqref{QP1} and \eqref{QP3},  the dominated convergence theorem gives
\beg{align*} \lim_{\vv\downarrow 0} \ff{f(\L_{X+\vv Y})- f(\L_{X}) }{\vv} &= \lim_{\vv\downarrow 0} \ff 1\vv\int_0^\vv \E \<D^Lf(\L_{X+sY}) (X+sY), Y\>\,\d s \\
&=\E \<D^Lf(\L_{X}) (X), Y\>.\end{align*}
 \end{proof}

\subsection{Partial  derivative in initial value}

  For any $T>0$, let $\C_T=C([0,T]\to\R^d)$ be the path space over $\R^d$ with time interval $[0,T],$ and let $X_0,\eta\in L^2(\OO\to\R^d,\F_0,\P)$. For any $\vv\ge 0$, let $(X_{t}^{\vv})_{t\ge 0}$ solve the SDE
\beq\label{E2} \d X_t^{\vv} = b_t(X_t^\vv, \L_{X_t^\vv}) \d t +\si_t(X_t^\vv, \L_{X_t^\vv})\d W_t,\ \ X_0^\vv= X_0+\vv \eta.\end{equation}
Obviously, $X_t=X_t^0$ solves \eqref{E1} with initial value $X_0$. Consider the following linear SDE for $v_t^\eta$ on $\R^d$:
\beq\label{DR2'} \beg{split} \d v_t^\eta & =\Big\{\nn_{v_t^\eta} b_t(\cdot, \L_{X_t})(X_t)+ \big(\E \< D^L b_t(y,\cdot)(\L_{X_t})(X_t), v_t^\eta\>\big)\big|_{y=X_t}\Big\}\d t\\
&+ \Big\{\nn_{v_t^\eta} \si_t(\cdot, \L_{X_t})(X_t)+ \big(\E \< D^L \si_t(y,\cdot)(\L_{X_t})(X_t), v_t^\eta\>\big)\big|_{y=X_t}\Big\}\d W_t,\ \
v_0^\eta=\eta.\end{split}\end{equation}

  The main result of this part is the following.

\beg{prp}\label{P2.1} Assume {\bf (H)}. Then for any $T>0$, the limit
 \beq\label{LD0} \nn_{\eta} X_t  := \lim_{\vv\downarrow 0} \ff{X_t^\vv-X_t }\vv,\ \ t\in [0,T]\end{equation} exists in $L^2(\OO\to \C_T,\P).$
 Moreover, $(v_t^\eta:=\nn_{\eta} X_t)_{t\in [0,T]}$ is the unique solution to the linear SDE $\eqref{DR2'}$.
 \end{prp}

To prove the existence of $\nn_\eta X_t$ in \eqref{LD0}, it suffices to show that when $\vv\downarrow 0$
\beq\label{PP0} \xi^\vv(t):= \ff{X_t^\vv-X_t}\vv,\ \ t\in [0,T]\end{equation}  is a Cauchy sequence in $L^2(\OO\to \C_T,\P),$  i.e.
\beq\label{CS} \lim_{\vv,\dd\downarrow 0} \E \bigg[\sup_{t\in [0,T]} |\xi^\vv(t)- \xi^\dd(t)|^2\bigg]=0.\end{equation}
To this end, we need the following two  lemmas.

\beg{lem}\label{L2.2} Assume {\bf (H)}. Then
\beq\label{LM1}  \sup_{\vv\in (0,1]} \E\bigg[\sup_{t\in [0,T]} |\xi^\vv(t)|^2\bigg]<\infty. \end{equation} 
Moreover, for any $t\in [0,T]$, $\{|\xi^\vv(t)|\}_{\vv\in (0,1]}$ is uniformly integrable in $L^2(\P)$. 
\end{lem}
\beg{proof} By {\bf (H)}, there exists a constant $C_1>0$ such that
\beg{align*} &\d |X_t^\vv-X_t|^2\\
 &=\big\{2\<b_t(X_t^\vv,\L_{X_t^\vv})- b_t(X_t,\L_{X_t}), X_t^\vv-X_t\> +\|\si_t(X_t^\vv,\L_{X_t^\vv})-\si_t(X_t,\L_{X_t})\|_{HS}^2\big\}\d t +\d M_t\\
&\le C_1\big\{|X_t^\vv-X_t|^2+\W_2(\L_{X_t^\vv}, \L_{X_t})^2\big\} \d t  +\d M_t,\end{align*}
where $$\d M_t:= 2 \Big\<X_t^\vv-X_t, (\si_t(X_t^\vv,\L_{X_t^\vv})- \si_t(X_t,\L_{X_t}))\d W_t\Big\>$$ satisfies
\beq\label{*DDO} \d\<M\>_t \le C_1^2 \big\{|X_t^\vv-X_t|^2+\W_2(\L_{X_t^\vv}, \L_{X_t})^2\big\}^2\d t.\end{equation}
Then by the BDG inequality, and noting that $\W_2(\L_\xi,\L_\eta)^2\le \E|\xi-\eta|^2$ for two random variables $\xi,\eta$, we may find out a constant  $C_2>0$ such that
\beq\label{*DDP}  \E\bigg[ \sup_{s\in [0,t]} |X_s^\vv-X_s|^2\bigg]  \le \vv^2|\eta|^2 + 2 C_1 \int_0^t \E|X_s^\vv-X_s|^2\d s + C_2 \E \ss{\<M\>_t}.\end{equation}
 Noting that $\W_2(\L_{X_s^\vv},\L_{X_s})^2\le \E|X_s^\vv-X_s|^2$, \eqref{*DDO} yields
\beg{align*} &C_2 \E\ss{\<M\>_t} \le C_1C_2 \E\bigg(\int_0^t \big\{|X_s^\vv-X_s|^2+\W_2(\L_{X_s^\vv}, \L_{X_s})^2\big\}^2\d s\bigg)^{\ff 1 2}\\
&\le C_1C_2 \E\bigg(\sup_{s\in [0,t]} \big\{|X_s^\vv-X_s|^2+\E|X_s^\vv- X_s|^2\big\} \int_0^t \big\{|X_s^\vv-X_s|^2+\E|X_s^\vv-X_s|^2\big\} \d s\bigg)^{\ff 1 2}\\
&\le \ff 1 2 \E \Big[\sup_{s\in [0,t]} |X_s^\vv-X_s|^2\Big] + \ff {C_3}2  \int_0^t \E|X_s^\vv-X_s|^2\,\d s \end{align*} for some constant $C_3>0$.
Combining this with \eqref{*DDP} and noting that due to \eqref{MM}  $$\E \Big[\sup_{s\in [0,t]} |X_s^\vv-X_s|^2\Big]<\infty,$$  we arrive at
$$\E\bigg[ \sup_{s\in [0,t]} |X_s^\vv-X_s|^2\bigg] \le 2 \vv^2|\eta|^2 + C_3 \int_0^t \E|X_s^\vv-X_s|^2\d s,\ \ t\in [0,T],\vv>0.$$
Therefore, Gronwall's inequality gives
$$\sup_{\vv\in (0,1]} \E \bigg[\sup_{t\in [0,T]} |\xi^\vv(t)|^2\bigg] =  \sup_{\vv\in (0,1]} \ff 1 {\vv^2} \E \bigg[\sup_{s\in [0,T]}  |X_s^\vv-X_s|^2\bigg]\le 2\e^{C_3T} \E|\eta|^2  <\infty,$$
which implies \eqref{LM1}. Moreover, there exists constant $C>0$ such that for any $\vv\in (0,1]$,
\beq\label{UNF} \d |\xi^\vv(t)|^2\le C(|\xi^\vv(t)|^2+1) \d t + \d M^\vv(t),\ \ t\ge 0\end{equation}
for some martingale $M^\vv(t)$ with $\d \<M^\vv\>(t)\le C\d t.$ 

To prove the uniform integrability of  $\{|\xi^\vv(t)|\}_{\vv\in (0,1]}$  in $L^2(\P)$, by \eqref{LM1} we find a constant $c_1>0$ such that 
$$\sup_{\vv\in (0,1]} \E\int_0^t |\xi^\vv(s)|^2 \d s \le c_1<\infty.$$
Then for any $N\in \mathbb N$ and $\vv\in (0,1]$, there exists $n\in [N, 2N]\cap\mathbb N$ such that 
\beq\label{U1} \E\int_0^t |\xi^\vv(s)|^2 1_{\{2^n\le |\xi^\vv(s)|^2\le 2^{n+1}\}} \d s \le \ff {c_1}N.\end{equation}
Take $h\in C_b^2([0,\infty)$ with $h'\ge 0, h|_{[0,1]}=0$ and $h|_{2,\infty)}=1$, and let  
$$h_n(r) = r h(2^{-n}r),\ \  r\ge 0.$$
By \eqref{UNF} and It\^o's formula, we find  a constant  $c_2>0$ such that 
\beg{align*} &\d h_n(|\xi^\vv(s)|^2) \\
&\le \big\{ C h_n(|\xi^\vv(s)|^2)+ C h(2^{-n} |\xi^\vv(s)|^2) \E|\xi^\vv(s)|^2 +  c_2 |\xi^\vv(s)|^2 1_{\{2^n\le |\xi^\vv(s)|^2\le 2^{n+1}\}}\big\}\d s + \d M^{\vv,n}(s)\end{align*}
for some martingale $M^{\vv,n}(s)$. Combining this with \eqref{U1} and noting that the FKG inequality implies
$$ \E[h(2^{-n} |\xi^\vv(s)|^2) ] \E[|\xi^\vv(s)|^2]\le \E h_n(|\xi^\vv(s)|^2),$$
we obtain
$$\E[h_n(|\xi^\vv(s)|^2)] \le \E [h_n(|\eta|^2)]  + \ff{c_1}N + 2C\int_0^s \E[h_n(|\xi^\vv(r)|^2)] \d r,\ \ s\in [0,t].$$ 
By Gronwall's lemma and the definition of $h_n$ for $n\in [N,2N]$, we arrive at
$$\E[|\xi^\vv(t)|^2 1_{\{|\xi^\vv(t)|^2\ge 4^N\}}] \le \E [h_n((|\xi^\vv(t)|^2)] \le \big\{\E [ |\eta|^21_{\{|\eta|^2\ge 2^N\}}]  +  c_1 N^{-1} \big\} \e^{2Ct},$$
which goes to zero uniformly in $\vv\in (0,1]$ as $N\to\infty$. 
Therefore,   $\{|\xi^\vv(t)|\}_{\vv\in (0,1]}$ is uniformly integrable in $L^2(\P)$. 
\end{proof}

For any differentiable (real, vector, or matrix valued) function $f$ on $\R^d\times \scr P_2(\R^d)$,  let
\beq\label{XII} \beg{split} \Xi^\vv_f(t)= &\ \ff{f(X_t^\vv,\L_{X_t^\vv})-f(X_t, \L_{X_t})}\vv -\nn_{\xi^\vv(t)} f(\cdot, \L_{X_t})(X_t) \\
&\quad -\big\{\E\<D^L f(y, \cdot)(\L_{X_t})(X_t),\xi^\vv(t)\>\big\}\big|_{y= X_t},\ \ t\in [0,T], \vv>0.\end{split}\end{equation}

\beg{lem}\label{L2.3} Assume {\bf (H)}. For any  $($real, vector, or matrix valued$)$ $C^{1,(1,0)}$-function $f$ on $\R^d\times \scr P_2(\R^d)$ with
\beq\label{BDD}K_f:= \sup_{(x,\mu)\in\R^d\times\scr P_2(\R^d)}\big(|\nn f(\cdot,\mu)(x)|^2+\|D^L f(x,\cdot)(\mu)\|_{L^2(\mu)}^2\big)<\infty,\end{equation}
there holds
\beq\label{BDD2} \big|\Zeta_f^\vv(t)\big|^2\le 4 K_f \big(\E |\xi^\vv(t)|^2 + |\xi^\vv(t)|^2\big) \ \ \text{and}\ \
 \lim_{\vv\downarrow 0} \E\big|\Zeta_f^\vv(t)\big|^2=0, \ \ \  t\in [0,T].\end{equation}
 \end{lem}
 \beg{proof} Let $X_t^{\vv}(s)= X_t+ s(X_t^\vv-X_t),\ s\in [0,1].$ By the chain rule and \eqref{XMM2},  we have
 \beg{align*} & \ff{f(X_t^\vv,\L_{X_t^\vv})-f(X_t, \L_{X_t})} \vv
 = \ff 1 \vv \int_0^1 \Big\{\ff{\d}{\d s} f\big(X_t^{\vv}(s), \L_{X_t^{\vv}(s)}\big)\Big\}\,\d s\\
 &= \int_0^1 \Big\{ \nn_{\xi^\vv(t)} f(\cdot, \L_{X_t^{\vv}(s)})(X_t^{\vv}(s)) + \big(\E\big\<D^L f(y,\cdot)(\L_{X_t^{\vv}(s)})(X_t^{\vv}(s)), \xi^\vv(t)\big\>\big)\big|_{y=X_t^{\vv}(s)}\Big\}\d s.\end{align*}
 Combining this with \eqref{BDD} we obtain
 \beq\label{PP}\beg{split}   \big|\Xi^\vv_f(t)\big|^2 \le &\ 2\int_0^1 \Big| \nn_{\xi^\vv(t)} \big\{f(\cdot, \L_{X_t^{\vv}(s)})(X_t^{\vv}(s))- f(\cdot, \L_{X_t})(X_t)\big\}\Big|^2\d s\\
 &+ 2 \int_0^1 \Big|  \big(\E\big\<D^L f(y,\cdot)(\L_{X_t^{\vv}(s)})(X_t^{\vv}(s)), \xi^\vv(t)\big\>\big)\big|_{y=X_t^{\vv}(s)}\\
  &\qquad \qquad\quad  -\big(\E\big\<D^L f(y,\cdot)(\L_{X_t})(X_t), \xi^\vv(t)\big\>\big)\big|_{y=X_t} \Big|^2\d s\\
  &\le 8 K_f (|\xi^\vv(t)|^2+\E|\xi^\vv(t)|^2\big).\end{split} \end{equation}
  So, the first inequality in \eqref{BDD2} holds.   Moreover, Lemma \ref{L2.2} implies
 $$\lim_{\vv\downarrow 0} \E \bigg[\sup_{s\in [0,1]} |X_t^\vv(s)-X_t|^2\bigg] \le \lim_{\vv\downarrow 0}\E  |X_t^\vv-X_t|^2=0.$$ Thus,  the $C^{1,(1,0)}$-property of $f$, Lemma \ref{L2.2} and the first inequality in \eqref{PP} yield that
 $\Xi^\vv_f(t)\to 0$ in probability as $\vv\to 0$. Combining this with the first inequality in \eqref{BDD2},  Lemma \ref{L2.2},  and using the dominated convergence theorem, we derive  $ \lim_{\vv\downarrow 0} \E\big|\Zeta_f^\vv(t)\big|^2=0.$ \end{proof}

\beg{proof}[Proof  of  Proposition \ref{P2.1}]   Let $(\Zeta_b^\vv(t), K_{b_t})$ and $(\Xi^\vv_\si(t), K_{\si_t})$ be defined as in \eqref{XII} and \eqref{BDD}  for $b_t$ and $\si_t$ replacing $f$ respectively.
By {\bf (H)}, there exists a constant $C_1>0$ such that
$$\sup_{t\in [0,T]} \big(K_{b_t} +K_{\si_t}\big)\le C_1<\infty.$$ Then Lemma \ref{L2.3} gives
\beq\label{PP2}\beg{split} & \big|\Xi^\vv_b(t)\big|^2+  \big|\Xi^\vv_\si(t)\big|^2\le 4C \big(|\xi^\vv(t)|^2 +\E |\xi^\vv(t)|^2\big),\\
&  \lim_{\vv\downarrow 0} \E  \big(\big|\Xi^\vv_b(t)\big|^2+  \big|\Xi^\vv_\si(t)\big|^2\big)=0,\ \ t\in [0,T]. \end{split}\end{equation}
By \eqref{E2}, \eqref{PP0}, and \eqref{XII}  for $b_t$ and $\si_t$ replacing $f,$  we have
\beg{align*} \xi^\vv(t) &= \int_0^t\Big\{\Xi^\vv_b(s) +\nn_{\xi^\vv(s)} b_s(\cdot, \L_{X_s})(X_s) + \big(\E\<D^L b_s(y, \cdot)(\L_{X_s}) (X_s), \xi^\vv(s)\> \big)\big|_{y=X_s}\Big\}\d s\\
& \quad +  \int_0^t\Big\<\Xi^\vv_\si(s) +\nn_{\xi^\vv(s)} \si_s(\cdot, \L_{X_s})(X_s) + \big(\E\<D^L \si_s(y, \cdot)(\L_{X_s}) (X_s), \xi^\vv(s)\> \big)\big|_{y=X_s},\ \d W_s\Big\>
  \end{align*}  for $t\in [0,T].$ So, for any $\vv,\dd\in (0,1],$  $\xi^{\vv,\dd}(t):= \xi^\vv(t)- \xi^\dd(t)$ satisfies
  \beg{align*} |\xi^{\vv,\dd}(t)|^2& \le 4 \int_0^t \big|\Zeta_b^\vv(s)- \Zeta_b^\dd(s)\big|^2 \d s+ 4 \bigg|\int_0^t \big\< \Zeta_\si^\vv(s)- \Zeta_\si^\dd(s), \d W_s\big\>\bigg|^2\\
  & + 4T \int_0^t \Big| \nn_{\xi^{\vv,\dd}(s)} b_s(\cdot, \L_{X_s})(X_s) + \big(\E\< D^L b_s(y, \cdot)(\L_{X_s})(X_s), \xi^{\vv,\dd}(s)\>\big)|_{y=X_s}\Big|^2\d s\\
  &+ 4\bigg| \int_0^t \Big\<\nn_{\xi^{\vv,\dd}(s)} \si_s(\cdot, \L_{X_s})(X_s) + \big(\E\< D^L \si_s(y, \cdot)(\L_{X_s})(X_s), \xi^{\vv,\dd}(s)\>\big)|_{y=X_s},\ \d W_s\Big\>\bigg|^2.
  \end{align*}
Combining this with {\bf (H)} and using the Burkholder-Dave-Gundy  inequality, we find out a constant $C_2>0$ such that
\beg{align*} \E\bigg[\sup_{s\in [0,t]}  \xi^{\vv,\dd}(s) \bigg] \le &\ C_2 \int_0^T \E\Big(\big|\Zeta_b^\vv(s)- \Zeta_b^\dd(s)\big|^2 + \big|\Zeta_\si^\vv(s)- \Zeta_\si^\dd(s)\big|^2 \Big)\d s \\
 &\qquad  + C_2 \int_0^t \E |\xi^{\vv,\dd}(s)|^2\ \d s,\ \ t\in [0,T].\end{align*}
Since Lemma \ref{L2.2} ensures that $\E\big[\sup_{s\in [0,t]}  \xi^\vv(s) \big]<\infty$, by Gronwall's lemma this yields
$$\E\bigg[\sup_{s\in [0,T]}  \xi^{\vv,\dd} (s) \bigg]\le C_2\e^{C_2 T} \int_0^T \E\Big(\big|\Zeta_b^\vv(s)- \Zeta_b^\dd(s)\big|^2 + \big|\Zeta_\si^\vv(s)- \Zeta_\si^\dd(s)\big|^2 \Big)\d s.$$
Combining this with \eqref{PP2} and Lemma \ref{L2.2},  and applying the dominated convergence theorem, we prove the first assertion in Proposition \ref{P2.1}.

Finally, by \eqref{E2}, \eqref{LD0}, \eqref{PP2} and \eqref{XII} for $b_t,\si_t$ replacing $f$,  we conclude that $v_t^\eta:= \nn_{\eta} X_t$ solves the SDE \eqref{DR2'}. Since this SDE is linear,    the uniqueness is trivial. Then the proof is finished.
\end{proof}

 \subsection{Malliavin derivative}

Consider the Cameron-Martin space
$$\H= \bigg\{h\in C([0,T]\to \R^d): h_0={\bf0}, h_t'\ \text{exists\ a.e.}\ t, \|h\|_\H^2:= \int_0^T | h_t'|^2\d t<\infty\bigg\}.$$
  Let $\eta\in L^2(\OO\to\R^d,\F_0,\P)$ with $\L_{\eta}=\mu,$  and let   $\mu_T$ be the distribution of $W_{[0,T]}:=\{W_t\}_{t\in [0,T]}$, which is a probability measure (i.e. Wiener measure) on the path space  $\C_T:=C([0,T]\to\R^d)$.
For  $F\in L^2(\R^d\times \C_T,\mu\times \mu_T)$, $F(\eta, W_{[0,T]})$  is called Malliavin differentiable along direction  $h\in \H$, if the directional derivative
$$D_h F(\eta, W_{[0,T]}):= \lim_{\vv\to 0} \ff{F(\eta, W_{[0,T]}+\vv h)-F(\eta, W_{[0,T]})}{\vv}$$ exists in $L^2(\OO,\P)$. If the map $\H\ni h\mapsto D_h F\in L^2(\Omega,\mu)$ is bounded,
then there exists a unique   $DF(\eta, W_{[0,T]})\in L^2(\OO\to\H,\P)$ such that $\<DF(\eta, W_{[0,T]}), h\>_\H= D_h F(\eta, W_{[0,T]})$ holds
in $L^2(\Omega,\P)$ for all $h\in\H$. In this case, we write $F(\eta, W_{[0,T]})\in \D(D)$ and call $DF(\eta, W_{[0,T]})$
the Malliavin gradient of $F(\eta, W_{[0,T]})$. It is well known that $(D,\D(D))$ is a closed linear operator from $L^2(\OO, \F_T, \P)$ to $L^2(\OO\to\H, \F_T,\P)$. The adjoint operator $(D^*,\D(D^*))$ of $(D,\D(D))$ is called Malliavin divergence.   For simplicity, in the sequel we   denote $F(\eta, W_{[0,T]})$ by $F$. Then we have the integration by parts formula
\beq\label{INT}\E\big(D_h F\big|\F_0\big)  =   \E\big(FD^*(h)\big|\F_0\big),\ \ \ F\in \D(D), h\in
\D(D^*).\end{equation}
It is well known that for adapted $h\in L^2(\OO\to \H,\P)$, one has $h\in \D(D^*)$ with
\beq\label{DD*} D^*(h)= \int_0^T \<h_t', \d W_t\>.\end{equation}  For more details and applications on Malliavin calculus one may refer to \cite{NN} and references therein.

To calculate the Malliavian derivative of $X_t$ with $\L_{X_0}=\mu\in \scr P_2(\R^d)$, we  write $X_t=F_t(W_\cdot)$ as a functional of the Brownian motion  $\{W_{s}\}_{s\in [0,t]}$. Then by definition, for an adapted $h\in L^2(\OO\to\H,\P)$,
$$D_h X_t =\lim_{\vv\downarrow  0}\ff{F_t(W_\cdot +\vv h_\cdot) -F_t(W_\cdot)}\vv,\ \ 0\le t\le T.$$
On the other hand, by the pathwise uniqueness of \eqref{E1},   see  for instances \cite{PECR19,SZ,W18},  
$X_t^{h,\vv}:=F_t(W_\cdot+\vv h_\cdot)$ solves the SDE
\beq\label{E3} \d X_t^{h,\vv}=  b_t(X_t^{h,\vv}, \L_{X_t})   \d t +\si_t(X_t^{h,\vv}, \L_{X_t})\d (W_t+\vv h_t),\ \ X_0^{h,\vv}= X_0,\end{equation}
which is well-posed due to {\bf (H)} and  $h_\cdot'\in L^2(\OO\times [0,T], \P\times \d t)$.
When $\si_t(x,\mu)$ does not depend $(x,\mu)$, this SDE reduces to a   random ODE for
$Y_t^{h,\vv}:= X_t^{h,\vv}-\si_tW_t$, which  is  well-posed  also for non-adapted $h$ like $h^\aa$  in Theorem   \ref{T4.2}.  The main result of this part is the following which is well known by   regarding \eqref{E1} as the classical SDE, since in \eqref{E3} the distribution  $\L_{X_t}$ does not depend on the variable $\vv$.

\beg{prp}\label{P2.4} Assume {\bf (H)}. Let  $h\in L^2(\OO\to\H,\P)$, which is adapted if $\si_t(x,\mu)$   depends on $x$ or $\mu$. Then the limit
 \beq\label{LD} D_h X_t := \lim_{\vv\downarrow 0} \ff{X_t^{h,\vv}-X_t}\vv ,\ \ t\in [0,T]\end{equation} exists in $L^2(\OO\to \C_T,\P).$
 Moreover, $(w_t^h:=D_h X_t)_{t\in [0,T]}$ is the unique solution to the   SDE
\beq\label{DR2}  \beg{split} \d w_t^h  =  & \Big\{\nn_{w_t^h} \si_t(\cdot, \L_{X_t})(X_t)\Big\}\d W_t\\
 &+ \Big\{\nn_{w_t^h} b_t(\cdot, \L_{X_t})(X_t)+\si_t(\cdot, \L_{X_t})(X_t)h_t'\Big\}\d t,\ \ w^h_0={\bf 0}.\end{split}  \end{equation}   \end{prp}


\section{Proofs of main results}

We first present an integration by parts  formula for $\nn_\eta X_T$ with $\eta\in L^2(\OO\to\R^d,\F_0,\P)$, then prove Theorem \ref{T3.1}, Corollary \ref{T3.2} and Theorem \ref{T4.2} respectively.
\subsection{An integration by parts    formula}

\beg{thm}\label{TNN} Assume {\bf (H)} and $\eqref{LLP}$. Then for any $f\in C^1_b(\R^d)$,  $\eta\in L^2(\OO\to\R^d,\P)$, and any  $0\le r<T$ and $g\in C^1([r,T])$ with $g_r=0$ and $g_T=1$,
\beq\label{BSM2'}\E\big(\<\nn f (X_T), \nn_\eta X_T\>\big|\F_r\big) =   \E\bigg(f(X_T)   \int_r^T \big\<\zeta_t^\eta,\ \d W_t\big\>\bigg|\F_r\bigg) \end{equation}  holds for
$$\zeta_t^\eta:=\si_t(X_t)^{-1} \Big\{g_t'v_t^\eta +\big(\E\<D^L b_t(y,\cdot)(\L_{X_t})(X_t), g_t v_t^\eta\>\big)\big|_{y= X_t}\Big\},\ \ t\in [0,T].$$\end{thm}

\beg{proof} Having Propositions \ref{P2.1} and  \ref{P2.4} in hands, the proof is more or less standard.
For $v_t^\eta$ solving $\eqref{DR2'}$, we take
\beq\label{LLP2} h_t=   \int_{t\land r}^t 1_{\{s\ge r\}} \zeta_s\,\d s,\ \ t\in [0,T].\end{equation}  By {\bf (H)}, \eqref{LLP},   and that $h\in L^2(\OO\to\H, \P)$  is adapted,     Proposition \ref{P2.4} applies.
Let $\tt v_t= g_t v_t^\eta$ for $t\in [r,T]$. Then  \eqref{DR2'}  and \eqref{LLP2} imply
\beg{align*} \d\tt v_t  & =\Big\{\nn_{\tt v_t} b_t(\cdot, \L_{X_t})(X_t)+ \big(\E \< D^L b_t(y,\cdot)(\L_{X_t})(X_t), \tt v_t\>\big)\big|_{y=X_t} +  g_t'  v_t^\eta \Big\}\d t\\
&\quad + \Big\{\nn_{\tt v_t} \si_t(\cdot, \L_{X_t})(X_t)\Big\}\d W_t\\
&= \Big\{\nn_{\tt v_t} b_t(\cdot, \L_{X_t})(X_t)
 +\si_t(X_t, \L_{X_t} )   h_t'\Big\}\d t
 +\Big\{\nn_{\tt v_t} \si_t(X_t)  \Big\}\d W_t,
\ \ t\ge r,\ \tt v_r={\bf 0}.\end{align*} So, $(\tt v_t)_{t\ge r}$  solves the SDE \eqref{DR2} with    $\tt v_r={\bf 0}$. On the other hand, by \eqref{LLP2} we have $h_t'=0$ for $t<r$, so that the solution to \eqref{DR2} with $w_0^h=0$ satisfies $w_r^h=0$. 
So,   the uniqueness of this SDE from time $r$ implies $\tt v_t= w_t^h$ for all $t\ge r$. Combining this with  Propositions \ref{P2.1} and \ref{P2.4},  we obtain
$$ \nn_\eta X_T= v_T^\eta= g_T v_T^\eta = \tt v_T =w_T^h= D_h X_T.$$
Thus, by the chain rule and the integration by parts formula \eqref{INT}, for any bounded $\F_r$-measurable $G\in \D(D)$, we have
\beg{align*} &  \E\big(G\<\nn f(X_T), \nn_\eta X_T\> \big)=
 \E\big(G\<\nn f(X_T), D_h X_T\>\ \big)= \E\big( G D_h f(X_T) \big) \\
 &= \E\big(D_h\{G f(X_T)\}- f(X_T) D_h G ) = \E(G f(X_T) D^*(h) \big),\end{align*} where in the last step we have used $D_hG=0$ since $G$ is $\F_r$-measurable but $h_t'=0$ for $t\le r$. Noting that the class of bounded $\F_r$-measurable $G\in\D(D)$ is dense in $L^2(\OO,\F_r,\P)$, this implies
 $$\E\big(\<\nn f(X_T), \nn_\eta X_T\>\big|\F_r \big)= \E(f(X_T) D^*(h) \big|\F_r\big).$$
 Combining this  with
$$D^*(h)= \int_r^T \< h_t',\d W_t\>= \int_r^T\big\<\zeta_t^\eta, \d W_t\big\>$$ due to \eqref{DD*} and \eqref{LLP2},
we prove  \eqref{BSM2'}.
\end{proof}

\subsection{Proof of Theorem \ref{T3.1}}

 Let $\mu\in \scr P_2(\R^d)$. We first establish \eqref{BSM1} for $f\in \B_b(\R^d)$, then construct   $\gg\in L^2(\R^d\to\R^d,\mu)$ such that
\beq\label{DLD} \lim_{\mu(|\phi|^2)\to 0} \ff{|(P_Tf)(\mu\circ({\rm Id}+\phi)^{-1})-(P_Tf)(\mu)- \mu(\<\phi,\gg\>)|}{\ss{\mu(|\phi|^2)}}=0,\end{equation}  which, by definition, implies that   $P_Tf$ is $L$-differentiable at $\mu$ with $D^L P_Tf(\mu) =\gg$.

 (a) Proof of \eqref{BSM1} for $f\in \B_b(\R^d)$.  When   $f\in C_b^1(\R^d)$, \eqref{BSM1} follows from \eqref{BSM2'} for $\eta=\phi(X_0)$. Below we extend the formula to $f\in \B_b(\R^d)$.
 For $s\in [0,1],$ let $X_t^{\phi,s}$ solve \eqref{E1} for $X_0^{\phi,s}=X_0+s\phi(X_0)$.
 We have  $\mu^{\phi,s}:= \L_{X_0^{\phi,s}}=\mu\circ({\rm Id}+s\phi)^{-1},$ and by   the definition of $\nn_\eta X_T$ for $\eta=\phi(X_0)$,
\beq\label{*DR} \beg{split} &(P_Tf)(\mu^{\phi,\vv})- (P_T f)(\mu)= \E[f(X_T^{\phi,\vv})-f(X_T)]= \int_0^\vv \ff{\d}{\d s} \E[f(X_T^{\phi,s})]\,\d s\\
&= \int_0^\vv   \E\<(\nn f)(X_T^{\phi,s}), \nn_{\phi(X_0)}X_T^{\phi,s}\>\,\d s,\ \ f\in C_b^1(\R^d).\end{split} \end{equation}
 Next, let  $(v_t^{\phi,s})_{t\in [0,T]}$ solve  \eqref{DR2'} for $\eta=\phi(X_0)$ and $X_t^{s}$ replacing $X_t$, i.e.
\beq\label{DR2''} \beg{split}  \d v_t^{\phi,s}  &=\Big\{\nn_{v_t^{\phi,s}} b_t(\cdot, \L_{X_t^{\phi,s}})(X_t^{\phi,s})+ \big(\E \< D^L b_t(y,\cdot)(\L_{X_t^{\phi,s}})(X_t^{\phi,s}), v_t^{\phi,s}\>\big)\big|_{y=X_t^{\phi,s}}\Big\}\d t\\
&+ \Big\{\nn_{v_t^{\phi,s}} \si_t(X_t^{\phi,s})\Big\}\d W_t,\ \ v_0^{\phi,s}=\phi(X_0).\end{split}\end{equation}   Let
$$  \zeta_t^{\phi,s}:= \si_t(X_t^{\phi,s})^{-1}\Big\{ g_t' v_t^{\phi,s}+ \big(\E \< D^L b_t(y,\cdot)(\L_{X_t^{\phi,s}})(X_t^{\phi,s}), g_t v_t^{\phi,s} \>\big)\big|_{y=X_t^{\phi,s}}\Big\},\ \ t\in [0,T].$$
Then \eqref{*DR} and   \eqref{BSM2'} imply
\beq\label{WPP01}  (P_Tf)(\mu^{\phi,\vv})- (P_T f)(\mu)
 =     \int_0^\vv  \E\bigg[f(X_T^{\phi,s})   \int_0^T \big\<  \zeta_t^{\phi,s},\ \d W_t\big\>\bigg]\,\d s,\ \ f\in C_b^1(\R^d), \end{equation}
By a standard approximation argument, we may extend  this formula to all $f\in \B_b(\R^d)$. Indeed, let
$$\nu_\vv(A)= \int_0^\vv \E\bigg[1_{A}(X_T^{\phi,s})   \int_0^T \big\<\zeta_t^{\phi,s},\ \d W_t\big\>\bigg]\,\d s,\ \ A\in \B(\R^d).$$
Then $\nu_\vv$ is a finite signed measure on $\R^d$  with
$$\int_{\R^d} f\d\nu_\vv= \int_0^\vv \E\bigg[f(X_T^{\phi,s})   \int_0^T \big\< \zeta_t^{\phi,s},\ \d W_t\big\>\bigg]\,\d s,\ \ f\in \B_b(\R^d).$$  So, \eqref{WPP01} is equivalent to
\beq\label{W001} \int_{\R^d}f\d P_T^*\mu^{\phi,\vv} - \int_{\R^d}f\d P_T^*\mu=\int_{\R^d} f\d\nu_\vv, \ \ f\in C_b^1(\R^d).\end{equation}
Since $\nu_{T,\vv}:=P_T^*\mu^{\phi,\vv}+P_T^*\mu+|\nu_\vv|$ is a finite measure on $\R^d$, $C_b^1(\R^d)$ is dense in $L^1(\R^d, \nu_{T,\vv})$. Hence,    \eqref{W001}  holds for all $f\in \B_b(\R^d)\subset L^1(\R^d, \nu_{T,\vv}).$ Consequently,
\eqref{WPP01} holds for all $f\in \B_b(\R^d).$ Thus,
\beq\label{WPP02}  \ff{ (P_Tf)(\mu^{\phi,\vv})- (P_T f)(\mu) }\vv
 =   \ff 1 \vv  \int_0^\vv  \E\bigg[f(X_T^{\phi,s})   \int_0^T \big\< \zeta_t^{\phi,s},\ \d W_t\big\>\bigg]\,\d s,\ \ f\in \B_b(\R^d). \end{equation}
 It is easy to see from {\bf (H)} that
$$\lim_{s\to 0} \sup_{t\in [0,T]} \E \big(|X_t^{\phi,s}-X_t|^2+ |v_t^{\phi,s}-v_t^\phi|^2\big)=0.$$
So,
\beq\label{*WP0} \lim_{\vv\downarrow 0} \ff 1 \vv \int_0^\vv \E \bigg| \int_0^T \big\<\zeta_t^{\phi,s}-\zeta_t^\phi,\d W_t\big\>\bigg| =0.\end{equation}
Combining this with \eqref{WPP02}, we see that \eqref{BSM1} for $f\in \B_b(\R^d)$ follows from
\beq\label{*WP} \lim_{\vv\downarrow 0} \E\bigg[ \{f(X_T^{\phi,\vv}) -f(X_T)\} \int_0^T \big\<\zeta_t^\phi,\d W_t\big\>\bigg] =0,\ \ f\in\B_b(\R^d).\end{equation}
To prove this equality,  we denote
$$I_r := \int_0^r \big\<\zeta_t^\phi,\d W_t\big\>,\ \ \ r\in (0,T).$$ Applying \eqref{BSM2'} with $ g_t:=\ff{t-r}{T-r}$ for $t\in [r,T]$ and using {\bf (H)}, we  may find out a constant $C(T,r)>0$ such that
\beg{align*} &\big| \E[ I_r \{f(X_T^{\phi,\vv}) -f(X_T)\}]\big| = \bigg|\E\bigg[I_r \int_0^\vv   \<\nn f(X_T^{\phi,s}), \nn_{\phi(X_0)} X_T^{\phi,s}\>\d s\bigg]\bigg| \\
&\le \E\bigg[|I_r|\cdot\bigg|\int_0^\vv \E\big(\<\nn f(X_T^{\phi,s}), \nn_{\phi(X_0)} X_T^{\phi,s}\>\big|\F_r\big)\d s\bigg|\bigg]\\
&\le \ff{C(T,r)}{T-r}\|f\|_\infty \int_0^\vv \E\bigg[|I_r|\bigg(\int_r^T \big| v_t^{\phi,s}\big|^2 \d t\bigg)^{\ff 1 2}\bigg]\d s,\ \ f\in C_b^1(\R^d).\end{align*}
By the  argument extending \eqref{WPP01} from $f\in C_b^1(\R^d)$ to $f\in \B_b(\R^d)$, we conclude from this that for any   $r\in (0,T)$,
\beg{align*} \lim_{\vv\downarrow 0} \sup_{\|f\|_\infty\le 1}  \big| \E[ I_r \{f(X_T^{\phi,\vv}) -f(X_T)\}]\big| =0.\end{align*}
   Therefore,
\beq\label{*WPN} \beg{split} &\limsup_{\vv\downarrow 0} \sup_{\|f\|_\infty\le 1}\bigg|\E\bigg[ \{f(X_T^{\phi,\vv}) -f(X_T)\} \int_0^T \big\<\zeta_t^\phi,\d W_t\big\>\bigg]\bigg|\\
&= \limsup_{\vv\downarrow 0} \sup_{\|f\|_\infty\le 1}\bigg|\E\bigg[ \{f(X_T^{\phi,\vv}) -f(X_T)\} \int_r^T \big\<\zeta_t^\phi,\d W_t\big\>\bigg]\bigg|\\
&\le 2     \bigg( \E \int_r^T |\zeta_t^\phi|^2 \d t \bigg)^{\ff 1 2},\ \ r\in (0,T).\end{split}\end{equation}
  By letting $r\uparrow T$ we prove \eqref{*WP}.

(b) For any $f\in \B_b(\R^d)$, we  intend to find out $\gg\in L^2(\R^d\to\R^d,\mu)$ such that
\beq\label{WPP1}  \E   \bigg[f(X_T)\int_0^T \big\<\zeta_t^\phi,\d W_t\big\>\bigg] =\mu(\<\phi,\gg\>),\ \ \phi\in L^2(\R^d\to\R^d,\mu). \end{equation}
When $f\in C_b(\R^d)$, in step (c) we will   deduce  from this and  \eqref{BSM1}  that  $\gg= D^LP_Tf(\mu)$.
To construct the desired $\gg$,   consider the SDE
$$ \d X_t^{\phi} = b_t(X_t^\phi, \L_{X_t^\phi}) \d t +\si_t(X_t^\phi)\d W_t,\ \ X_0^\phi= X_0+ \phi(X_0),$$
and let $v_t^\phi$ solve  \eqref{DR0}. Since  \eqref{DR0} is a linear equation for $v_t^\phi$ with initial value $\phi(X_0)\in L^2(\OO\to\R^d,\F_0,\P)$,
the functional
$$L^2(\R^d\to\R^d,\mu)\ni \phi\mapsto L\phi:=   \E   \bigg[f(X_T)\int_0^T \big\<\zeta_t^\phi,\d W_t\big\>\bigg]$$ is linear, and
 by {\bf (H)} and \eqref{LLP},  there exists a constant $C(T)>0$ such that
$$ |L\phi|^2\le   C(T)\, \E |\phi(X_0)|^2= C(T)\,\mu(|\phi|^2),\ \ \phi\in L^2(\R^d\to\R^d,\mu).$$  Then $L$ is a bounded linear functional on the Hilbert space $L^2(\R^d\to\R^d,\mu)$. By Riesz's representation theorem,    there exists a unique
$\gg\in L^2(\R^d\to \R^d,\mu)$ such that
$$L\phi= \mu(\<\gg,\phi\>),\ \ \phi\in L^2(\R^d\to\R^d,\mu).$$   Therefore,   \eqref{WPP1} holds.

(c) Now, for $f\in \B_b(\R^d)$, we intend to verify \eqref{DLD} for $\gg$ in  \eqref{WPP1}, so that $P_Tf$ is $L$-differentiable with $D^L(P_Tf)(\mu)= \gg$.
By \eqref{WPP02} for $\vv=1$, we have
\beq\label{WPP}  (P_Tf)(\mu^1)- (P_T f)(\mu)
 = \int_0^1  \E\bigg[f(X_T^{\phi,s})   \int_0^T \big\<\zeta_t^{\phi,s},\ \d W_t\big\>\bigg],\ \ f\in \B_b(\R^d).  \end{equation}
For $\R^d$ random variables $X,v$, let
$$N_t(X,v)= \si_t(X)^{-1}\Big\{ g_t' v+ \big(\E \< D^L b_t(y,\cdot)(\L_{X})(X), g_t v\>\big)\big|_{y=X}\Big\},\ \ t\in [0,T].$$
Then $\zeta_t^{\phi,s}= N_t(X_t^{\phi,s}, v^{\phi,s})$ and $\zeta_t^{\phi}= N_t(X_t, v^{\phi}).$
Combining this with \eqref{WPP1} and \eqref{WPP},  and noting that $\mu^1= \mu\circ({\rm Id}+\phi)^{-1})$,  we arrive at
\beq\label{WPP0} \ff{|(P_Tf)(\mu\circ({\rm Id}+\phi)^{-1}))-(P_Tf)(\mu)- \mu(\<\phi,\gg\>)|}{\ss{\mu(|\phi|^2)}}\le \ \vv_1(\phi)+\vv_2(\phi) +\vv_3(\phi), \end{equation}  where
\beg{align*} &\vv_1(\phi):= \ff {1}{\ss{\mu(|\phi|^2)}}  \int_0^1 \E \bigg|\big(f(X_T^{\phi,s})-f(X_T)\big)\int_0^T \<  \zeta_t^{\phi,s}, \d W_t\>\bigg|\d s,\\
&\vv_2(\phi):= \ff {\|f\|_\infty}{\ss{\mu(|\phi|^2)}}  \int_0^1 \E \bigg|   \int_0^T \big\<N_t(X_t^{\phi,s}, v^{\phi}) - N_t(X_t, v^{\phi}) , \d W_t\big\>\bigg|\d s,\\
&\vv_3(\phi):= \ff {\|f\|_\infty}{\ss{\mu(|\phi|^2)}}  \int_0^1 \E \bigg|   \int_0^T
\big\<N_t(X_t^{\phi,s}, v^{\phi,s}) -N_t(X_t^{\phi,s}, v^{\phi}) , \d W_t\big\>\bigg|\d s.
 \end{align*}

  It is easy to deduce from {\bf (H)} that for any $p\ge 2$ there exists a constant $c(p)>0$ such that
\beq\label{WPP3} \sup_{t\in [0,T],s\in [0,1]}\E\big(|X_t^{\phi,s}-X_t|^p+|v_t^{\phi,s}|^p\big|\F_0\big) \le c(p)  |\phi(X_0)|^p.  \end{equation} Combining this with
the   continuity  of   $\si_t(x)$ in $x$ uniformly in $t\in [0,T]$, we conclude  that
\beq\label{WWP4} \lim_{\mu(|\phi|^2)\to 0}  \vv_2(\phi)=0.\end{equation}  Next, by the argument deducing \eqref{BSM1} from \eqref{WPP02}, it is easy to see that \eqref{WPP3} implies
\beq\label{WWP4'} \lim_{\mu(|\phi|^2)\to 0}  \vv_1(\phi)=0.\end{equation}
 Moreover, by the SDEs for $v_t^{\phi,s}$ and $v_t^\phi$ we have
\beq\label{SSS} \d (v_t^{\phi,s}-v_t^\phi)= \big\{ A_t (v_t^{\phi,s}-v_t^\phi) + \tt A_t v_t^{\phi,s}\big\}\d t+ \big\{B_t (v_t^{\phi,s}-v_t^\phi)+ \tt B_t v_t^\phi\big\}\d W_t,\end{equation}  where for a square integrable random variable $v$ on $\R^d$,
\beg{align*} &A_t v:= \nn_{v} b_t(\cdot, \L_{X_t})(X_t)+ \big(\E \< D^L b_t(y,\cdot)(\L_{X_t})(X_t), v\>\big)\big|_{y=X_t},\\
&\tt A_t v:= \nn_{v}  b_t(\cdot, \L_{X_t^{\phi,s}})(X_t^{\phi,s})
  +\big(\E \< D^L b_t(y,\cdot)(\L_{X_t^{\phi,s}})(X_t^{\phi,s}), v\>\big)\big|_{y=X_t^{\phi,s}}\\
  &\qquad \qquad  - \nn_v b_t(\cdot, \L_{X_t})(X_t)- \big(\E \< D^L b_t(y,\cdot)(\L_{X_t})(X_t), v\>\big)\big|_{y=X_t},\\
& B_t v:=  \nn_{v} \si_t(X_t), \ \
 \tt B_tv:=  \nn_{v}  \si_t(X_t^{\phi,s})-   \nn_v\si_t(X_t).\end{align*}
Combining this with \eqref{WPP3} and {\bf (H)},  we find a constant $c_1>0$ such that  
\beq\label{W005}  \vv_t:=  \|\tt A_t\|  +\|\tt B_t\| \le c_1,\ \ \lim_{\mu(|\phi|^2)\to 0}  \vv_t =0, \ \ t\in [0,T], s\in [0,1],\end{equation}
and  
$$g_t:=  \E(v_t^{\phi,s}-v_t^\phi||\F_0),\ \ t\in [0,T]$$ satisfies 
\beq\label{SSS2}  g_t \le c_1  \int_0^t (g_r+\E g_r) \d r +  c_1 \bigg(\int_0^t \E(\vv_r^2 (|v_r^\phi|^2 + |v_r^{\phi,s}|^2)|\F_0) d r\bigg)^{\ff 1 2},\ \ t\in [0,T].\end{equation}
On the other hand, by \eqref{WPP3} for $p=4$, we find a constant $c_2>0$ such that 
\beq\label{SSS5} \E(\vv_s^2 (|v_r^\phi|^2 + |v_r^{\phi,s}|^2)|\F_0) \le \{\E(\vv_r^4|\F_0)\}^{\ff 1 2} |\phi(X_0)|^2,\end{equation} 
so that \eqref{SSS2} yields
\beq\label{SSS4} g_t \le c_1  \int_0^t (g_r+\E g_r) \d r +  c_3  \{\E(\vv_r^4|\F_0)\}^{\ff 1 4} |\phi(X_0)|,\ \ t\in [0,T]\end{equation} 
for some constant $c_3>0$. Taking expectation in both sides and applying  Gronwall's lemma, we arrive at
\beq\label{SSS3} \vv(\phi):=\sup_{t\in [0,T], s\in [0,1]}  \ff 1 {\|\phi\|_{L^2(\mu)}} \E |v_t^{\phi,s}-v_t^\phi|\to 0\ \text{as}\ \|\phi\|_{L^2(\mu)}\downarrow.\end{equation} 
By the same argument leading to \eqref{SSS4}, we find a constant $c_4>0$ such that
\beg{align*}& \E(|v_t^{\phi,s}-v_t^\phi|^2|\F_0)\\
&\le c_4 \int_0^t \big\{\E(|v_r^{\phi,s}-v_r^\phi|^2|\F_0)+ (\E|v_r^{\phi,s}-v_r^\phi|)^2+ \{\E(\vv_r^4|\F_0)\}^{\ff 1 2} |\phi(X_0)|^2\big\}\d r,\ \ t\in [0,T].\end{align*}  
Combining this with   \eqref{SSS} we find a constant $c_5>0$ such that 
$$\sup_{t\in [0,T]} \E(v_t^{\phi,s}-v_t^\phi|^2|\F_0)\le c_5  \big(\vv(\phi)^2\|\phi\|_{L^2(\mu)}^2 + \dd(\phi)^2 |\phi(x_0)|^2\big),$$
 where according to \eqref{W005},
\beq\label{SSS6}  \dd(\phi) := \bigg(\int_0^T  \{\E(|\vv_r^4|\F_0)\} \d r \bigg)^{\ff 1 4} \le c_1 T,\ \  \lim_{\|\phi\|_{L^2(\mu)}\to 0} \dd(\phi) =0.\end{equation} 
So, by  the definition of $\vv_3(\phi)$, {\bf (H)},      we   find   constants
 $c_5, c_6, c_7>0$ depending on $\|f\|_\infty$ and $T$  such that \eqref{SSS5}, \eqref{SSS4},  \eqref{SSS3} and \eqref{SSS6} imply 
 \beg{align*} &  \vv_3(\phi)  
 \le    \ff{c_5}{\|\phi\|_{L^2(\mu)}} \int_0^1 \E\bigg(\int_0^T \E\big(|N_t(X_t^{\phi,s}, v_t^{\phi,s})- N_t(X_t^{\phi,s}, v_t^\phi)|^2\big|\F_0\big)  \d t\bigg)^{\ff 1 2}\d s\\
&\le    \ff{c_6}{\|\phi\|_{L^2(\mu)}} \int_0^1 \E\bigg(\int_0^T\big\{ \E(|v_t^{\phi,s}-v_t^\phi|^2 +\vv_t^2(|v_t^{\phi,s}|^2+|v_t^\phi|^2)|\F_0)+(\E |v_t^{\phi,s}-v_t^\phi| )^2\big\}\d t\bigg)^{\ff 1 2}\d s \\
&\le c_7\Big(\vv(\phi)+ \ff{ \E [\dd(\phi)|\phi(X_0)|]}{\|\phi\|_{L^2(\mu)}} \Big)\to  0 \ \text{as}\  \|\phi\|_{L^2(\mu)}\downarrow 0.
  \end{align*}    This,  together with \eqref{WPP0}, \eqref{WWP4} and \eqref{WWP4'}, implies \eqref{DLD}. Therefore,    $P_Tf$ is $L$-differentiable at $\mu$ with $D^L(P_Tf)(\mu)=\gg$.

  (d) Finally, \eqref{BSM1} and \eqref{WPP02} imply
  \beg{align*} &\Big| \ff{P_T^*\mu\circ({\rm Id}+\vv\phi)^{-1}- P_T^*\mu}\vv (f)- (\psi P_T^*\mu)(f)\Big|\\
  &= \bigg|\ff{(P_Tf)(\mu^{\phi,\vv}) - (P_Tf)(\mu)}{\vv}- \E\bigg[f(X_T)\int_0^T\<\zeta_t^\phi, \d W_t\>\bigg]\bigg|\\
  &\le \ff {\|f\|_\infty}  \vv \int_0^\vv  \E\bigg|\int_0^T\<\zeta_t^{\phi,s}-\zeta_t^\phi, \d W_t\> \bigg|\d s\\
  &\quad + \ff 1 {\vv}\bigg|\E\bigg[ \{f(X_T^{\phi,\vv}) -f(X_T)\} \int_0^T \big\<\zeta_t^\phi,\d W_t\big\>\bigg]\bigg|\d s.\end{align*} Combining this with \eqref{*WP0} and \eqref{*WP}  we prove \eqref{*D*}.

\subsection{Proof of Corollary \ref{T3.2}}

 \beg{proof}[Proof of $(1)$]   By {\bf (H)} and \eqref{DR0}, there exists a martingale $M_t$ such that
  \beq\label{**1} \d |v_t^\phi|^2 \le 4K_t |v_t^\phi|(|v_t^\phi|+\E |v_t^\phi|)\d t  +\d M_t ,\ \ |v_0^\phi|^2= |\phi(X_0)|^2,\end{equation} where $K(t)$ is increasing in $t\ge 0$. Then
  $$  \E |v_t^\phi|^2\le \E |\phi(X_0)|^2 + 4K_t \int_0^t \big\{\E|v_s^\phi|^2+ (\E|v_s^\phi|)^2\big\}\d s\le \mu(|\phi|^2) + 8K_t \int_0^t \E|v_s^\phi|^2\d s.$$
  By Gronwall's inequality this implies
\beq\label{PPG} \E |v_t^\phi|^2 \le    \e^{8K_tt} \mu(|\phi|^2),\ \   t\in [0,T].\end{equation}   Next, since $\E \int_0^T \big\< \xi_t^\phi,\ \d W_t\big\>=0,$
\eqref{BSM1}  is equivalent to
$$D_\phi^L (P_T f)(\mu) = \E\bigg[\big\{f(X_T)  - P_T f(\mu)\big\} \int_0^T \big\<\zeta_t^\phi,\ \d W_t\big\>\bigg].$$
Combining this with \eqref{PPG}
and using Jensen's inequality,  when $\mu(|\phi|^2)\le 1$ we have
 \beg{align*} &|D_\phi^L(P_Tf)(\mu)|^2 \le \big\{(P_Tf^2)(\mu) -(P_Tf(\mu))^2\big\}\int_0^T \E\big|\zeta_t^\phi\big|^2\d t\\
& \le \big\{(P_Tf^2)(\mu) -(P_Tf(\mu))^2\big\} \int_0^T  \big(| g_t' | +K(t)|g_t|\big)^2\ll_t^2  \e^{8 t K_t}\d t \end{align*}for any $g\in C^1([0,T])$ with $g_0=0$ and $g_T=1$.
Taking
$g_t= \ff t T,\  t\in [0,T],$
we prove the estimate \eqref{GRD}.\end{proof}

\beg{proof}[Proof of $(2)$]
  Let $f\in \B_b(\R^d)$ with $\|f\|_\infty\le 1$. By Theorem \ref{T3.1}, $P_Tf$ is $L$-differentiable. Moreover, by Theorem \ref{TNN}, $P_Tf$ is Lipschitz continuous on $\scr P_2(\R^d)$. Indeed, for any $\mu_1,\mu_2\in \scr P_2(\R^d)$, let $X_1,X_2\in L^2(\OO\to\R^d,\F_0,\P)$ such that $\L_{X_i}=\mu_i, 1\le i\le 2,$ and $\E|X_1-X_2|^2=\W_2(\mu_1,\mu_2)^2$. Let $X_t^s$ be the solution to \eqref{E1} with $X_0=X_1 +s(X_2-X_1), s\in [0,1].$ Then Theorem \ref{TNN} implies
 \beg{align*} & |P_Tf(\mu_1)-P_Tf(\mu_2)|^2 = |\E f(X_T^0)-\E f(X_T^1)|^2=\bigg|\int_0^1 \ff{\d}{\d s} \E f(X_T^s)\,\d s\bigg|^2\\
 &=\bigg|\int_0^1 \E \<\nn f(X_T^s), \nn_{X_2-X_1}X_T^s\>\d s\bigg|^2
  \le c \E |X_2-X_1|^2= c \W_2(\mu_1,\mu_2)^2 \end{align*} for some constant $c>0$.

 To apply  Proposition \ref{P01},  we take $\{\mu_n,\nu_n\}_{n\ge 1}\subset \scr P_2(\R^d)$ which have compact supports and are absolutely continuous with respect to the Lebesgue measure, such that
 \beq\label{FY6}\lim_{n\to\infty} \big\{\W_2(\mu,\mu_n)+\W_2(\nu,\nu_n)\big\}=0.\end{equation}
 According to \cite{BRE}, see also  \cite[Theorem 5.8]{Card}, for any $n\ge 1$ there exists a unique map $\phi_n\in L^2(\R^d\to\R^d,\mu)$ such that
\beq\label{FY7} \nu_n=\mu_n\circ ({\rm Id}+ \phi_n)^{-1},\ \   \W_2(\mu_n,\nu_n)^2= \mu_n (|\phi_n|^2).\end{equation}
   Let $X_n\in L^2(\OO\to\R^d,\F_0,\P)$ such that $\L_{X_n}=\mu_n$. By Proposition \ref{P01},
   \eqref{GRD} and \eqref{FY7}, we obtain \beg{align*} & |(P_Tf)(\mu_n)-(P_Tf)(\nu_n)|^2=\bigg|\int_0^1\ff{\d}{\d s} (P_T f)(\L_{ X_n+s \phi_n(X_n)})\,\d s\bigg|^2\\
 &=\bigg|\int_0^1 \E\big\<D^L (P_T f)(\L_{X_n+s\phi_n(X_n)})(X_n+s\phi_n(X_n)), \phi_n(X_n)\big\> \,\d s\bigg|^2\\
 & \le \ff{ \|f\|_\infty^2\mu_n(|\phi_n|^2)}{\int_0^T \ll_t^{-2} \e^{-8tK_t}\d t}  = \ff{\|f\|_\infty^2 \W_2(\mu_n,\nu_n)^2}{\int_0^T \ll_t^{-2} \e^{-8tK_t}\d t}.\end{align*}
 By the continuity of $P_Tf$ and \eqref{FY6}, by  letting  $n\to\infty$  we prove
$$ |(P_Tf)(\mu)-(P_Tf)(\nu)|^2\le \ff{ \W_2(\mu,\nu)^2}{\int_0^T \ll_t^{-2} \e^{-8tK_t}\d t},\ \
 \mu,\nu\in \scr P_2(\R^d),\ \ f\in \B_b(\R^d), \|f\|_\infty\le 1.$$ Therefore, \eqref{GRD3} and \eqref{GRD2} hold.
 \end{proof}

\subsection{Proof of Theorem \ref{T4.2}} Let $T>r\ge 0, \mu\in\scr P_2(\R^{m+d})$ and let $X_t$ solve \eqref{E5} with $\L_{X_0}=\mu$.
To realize the procedure in the proof of Theorem \ref{T3.1} for the present degenerate setting, we first extend Theorem \ref{TNN} using $D^*(h_{r,\cdot}^\aa)$ to replace $\int_r^T \<\zeta_t^\eta,\d W_t\>$, where
for a $C^1([r,T]\to \R^{m+d})$-valued random variable $\aa_\cdot=(\aa^{(1)}_\cdot,\aa^{(2)}_\cdot)$, let $(h_{r,t}^\aa, w_{r,t}^\aa)_{t\in [r,T]}$ be the unique solution to the random ODEs
\beq\label{B00*}\beg{split}& \ff{\d h_{r,t}^\aa}{\d t}  =  \si_t^{-1} \Big\{\nn_{\aa_t} b_t^{(2)}(X_t,\L_{X_t})-(\aa_t^{(2)})'\\
 &\qquad\qquad + \big(\E\<D^L b_t^{(2)}(y,\cdot)(\L_{X_t})(X_t), \aa_t+w_{r,t}^\aa\>\big)\big|_{y=X_t}\Big\},\\
& \ff{\d w_{r,t}^\aa}{\d t}= \nn_{w_{r,t}^\aa} b_t (\cdot,\L_{X_t})(X_t) + ({\bf 0}, \si_t (h_{r,t}^\aa)'),\ \ \ \ h_{r,r}^\aa=0, w_{r,r}^\aa=0. \end{split}\end{equation}

\beg{thm} \label{T4.1} Assume ${\bf (H1)}$.  Let $T>r\ge 0$, $\eta \in L^2(\OO\to\R^{m+d},\F_0,\P)$, and let $X_t$ solve $\eqref{E5}$ with $\L_{X_0}=\mu\in \scr P_2(\R^{m+d})$.
If there exists a $C^1([r,T]\to \R^{m+d})$-valued random variable $\aa_\cdot=(\aa^{(1)}_\cdot,\aa^{(2)}_\cdot)$ such that $\aa_r =\nn_\eta X_r, \aa_T={\bf 0},$
\beq\label{AA} (\aa_t^{(1)})'= \nn_{\aa_t} b_t^{(1)}(X_t),\ \ t\in [r,T],\end{equation} and $h_{r,\cdot}^\aa\in \D(D^*),$
then for any $f\in C_b^1(\R^{m+d}),$
\beq\label{MLL} \E\big(\<\nn f(X_T), \nn_\eta X_T\>\big|\F_r\big) =   \E\big(f(X_T)\, D^*(h_{r,\cdot}^\aa)\big|\F_r\big).\end{equation}
\end{thm}

\beg{proof} Letting $w_{t}= w_{r,t}^\aa1_{\{t>r\}}$,  Proposition \ref{P2.4} implies that  $w_{t} = D_{h^\aa_{r,\cdot}} X_t,\ t\in [0,T].$ By \eqref{B00*}, we have
$$w_t= \int_{t\land r}^t \Big\{\nn_{w_s}b_s(\cdot,\L_{X_s})(X_s) + \Big({\bf 0}, \si_s  (h^\aa_{r,s})'  \Big\}\d s,\ \ t\in [0,T].$$
Extending $\aa_t$ with $\aa_t:= \nn_\eta X_t$ for $t\in [0,r)$,  and letting  $v_t= w_t+\aa_t$  for any $t\in [0,T]$,  we obtain
\beq\label{WR1} \beg{split} & v_t = \aa_t + \int_{t\land  r}^t \Big\{\nn_{v_s}b_s(\cdot,\L_{X_s})(X_s) + \Big({\bf 0},   \big(\E\<D^L b_s^{(2)}(y,\cdot)(\L_{X_s})(X_s), v_s\>\big)\big|_{y=X_s}\Big)  \\
&    +  ({\bf 0}, \si_s(h^\aa_s)' - \big(\E\<D^L b_s^{(2)}(y,\cdot)(\L_{X_s})(X_s), w_s+\aa_s\>\big)\big|_{y=X_s})- \nn_{\aa_s}b_s(\cdot,\L_{X_s})(X_s)  \Big\}\d s.\end{split}\end{equation}
By \eqref{AA},
$$\int_{t\land  r}^t \nn_{\aa_s}b_s^{(1)}(\cdot,\L_{X_s})(X_s)\,\d s  = 1_{\{t> r\}}\big(\aa_{t}^{(1)}-\nn_{\eta} X_r^{(1)}\big),$$
while the definition of $h_{r,s}^\aa$ implies
\beg{align*}&\int_{t\land  r}^t\Big\{\si_s(h^\aa_s)'- \big(\E\<D^L b_s^{(2)}(y,\cdot)(\L_{X_s})(X_s), w_s+\aa_s\>\big)\big|_{y=X_s}- \nn_{\aa_s}b_s^{(2)}(\cdot,\L_{X_s})(X_s)\Big\}\d s\\
&=-\int_{t\land  r}^t   (\aa_s^{(2)})'\d s=   1_{\{t> r\}}\big(\nn_\eta X_r^{(2)}- \aa_t^{(2)}\big).\end{align*}
Combining these with \eqref{WR1} and Proposition \ref{P2.1} leads to
\beg{align*} v_t &= \nn_\eta X_r  + \int_{t\land r}^t \Big\{\nn_{v_s}b_s(\cdot,\L_{X_s})(X_s) + \Big({\bf 0},   \big(\E\<D^L b_s^{(2)}(y,\cdot)(\L_{X_s})(X_s), v_s\>\big)\big|_{y=X_s}\Big) \Big\}\d s\\
&=  \eta   + \int_{0}^t \Big\{\nn_{v_s}b_s(\cdot,\L_{X_s})(X_s) + \Big({\bf 0},   \big(\E\<D^L b_s^{(2)}(y,\cdot)(\L_{X_s})(X_s), v_s\>\big)\big|_{y=X_s}\Big) \Big\}\d s ,\ \ t\in [0,T].\end{align*}
That is,   $v_t$ solves \eqref{DR2'} so that  by Proposition \ref{P2.1} we obtain $v_t:=w_t+\aa_t= \nn_\eta X_t.$ Since $\aa_T=0$, this implies
$D_{h_{r,\cdot}^\aa} X_T= \nn_\eta X_T$. Thus,   for any bounded $\F_r$-measurable $G\in \D(D)$,
\beq\label{GGF} \beg{split} &\E\big[ G \<\nn f(X_T), \nn_\eta X_T\>\big]= \E \big[G D_{h_{r,\cdot}^\aa} f(X_T)\big]\\
& = \E \big[D_{h_{r,\cdot}^\aa} \{Gf(X_T)\}- f(X_T)D_{h_{r,\cdot}^\aa} G \big]= \E \big[Gf(X_T)D^*(h_{r,\cdot}^\aa) \big],\end{split}\end{equation}
where in the last step we have used the integration by parts formula  \eqref{INT} and $D_{h_{r,\cdot}^{\aa}}G=0$ since $G$ is $\F_r$-measurable but
$$  D_{h_{r,\cdot}^\aa} G= \int_0^T (h_{r,\cdot}^\aa)' (s) \cdot \{(DG)_\cdot\}'(s)\d s=0,$$
  $(h_{r,\cdot}^\aa)' (s)=0$ for $s\le r$.   Noting that the class  of bounded $\F_r$-measurable functions $G\in \D(D)$ is dense in $L^2(\OO,\F_r,\P)$, \eqref{GGF} implies  \eqref{MLL}. \end{proof}

\beg{proof}[Proof of Theorem \ref{T4.2}]  With Theorem \ref{T4.1} in hands, the proof is completely similar to that of Theorem \ref{T3.1}.
Let $$v_t^\phi=((v_t^\phi)^{(1)}, (v_t^\phi)^{(2)})= (\nn_{\phi(X_0)}X_t^{(1)}, \nn_{\phi(X_0)}X_t^{(2)})= \nn_{\phi(X_0)}X_t,\ \ t\in [0,T].$$   For any $0\le r<T$,  let
\beg{equation}\label{aa'}\beg{split} \aa_{r,t}^{(2)}  = &\ff{T-t} {T-r} (v_t^\phi)^{(2)} -\ff{(t-r)(T-t) B_t^*K_{T,t}^*}{\int_0^T \theta_s^2 \d s} \int_t^T  \theta_s^2 Q_s^{-1} K_{T,r}(v_t^\phi)^{(1)}\d s\\
&-(t-r)(T-t) B_t^* K_{T,t}^*Q_T^{-1}\int_0^T\ff{T-s}TK_{T,s}\nn^{(2)}b^{(1)}_s(X_s) \phi^{(2)} (X_0)\d s,\ \ t\in [r,T],
\end{split}\end{equation} and
\beg{equation}\label{B0'} \aa^{(1)}_{r,t}= K_{t,r}(v_t^\phi)^{(1)}+\int_r^tK_{t,s}\nn^{(2)}_{\aa_s^{(2)}} b_s^{(1)}(X_s(x))\,\d s, \ \ t\in [r,T]. \end{equation}   Then $\aa_{r,\cdot}:=(\aa_{r,t}^{(1)},\aa_{r,t}^{(2)})$ satisfies
$$\aa_{r,r}= \nn_{\phi(X_0)}X_r,\ \ \aa_{r,T}=0,$$ and by \eqref{Eq1} and Duhamel's formula, \eqref{B0'} implies
$$(\aa_{r,\cdot}^{(1)})'(t)= \nn_{\aa_{r,t}} b_t^{(1)}(X_t),\ \ t\in [r,T].$$ Moreover, let $h_{r,\cdot}^{\aa_{r,\cdot}}$ be defined in \eqref{B00*} for $\aa_{r,\cdot}$ replacing $\aa$.
  Noting that {\bf (H1)} and {\bf (H2)} imply \cite[(H)]{WZ13} for $l_1=l_2=0$,  the proof   of \cite[Theorem 1.1]{WZ13} with $\phi(s):=(s-r)(T-s)$ for $s\in [r,T]$ ensures that $h_{r,\cdot}^{\aa_{r,\cdot}}\in\D(D^*)$ with $D^* (h_{r,\cdot}^{\aa_{r,\cdot}})\in L^p(\P)$ for all $p\in (1,\infty).$ So, by Theorem \ref{T4.2} with $\eta=\phi(X_0)$ we obtain
\beq\label{INTN}  \E (\<\nn f(X_T), \nn_{\phi(X_0)} X_T\>|\F_r) = \E(f(X_T) D^*(h_{r,\cdot}^{\aa_{r,\cdot}})|\F_r),\ \ f\in C_b^1(\R^d), r\in [0,T).\end{equation}
In particular, taking $r=0$ we obtain $D^*(h) \in L^p(\P)$ for all $p\in (1,\infty)$ and
\beq\label{INTN2}  D_\phi^LP_Tf(\mu)= \E (\<\nn f(X_T), \nn_{\phi(X_0)} X_T\>) = \E(f(X_T) D^*(h^{\aa})|\F_r),\ \ f\in C_b^1(\R^d).\end{equation}
Basing on these two formulas,    by repeating the proof of Theorem \ref{T3.1} with    $I_r:=\E(D^*(h^\aa)|\F_r)$, we prove
  \eqref{BSMN} and the $L$-differentiability of $P_Tf$ for $f\in \B_b(\R^{m+d})$. Finally, the estimates \eqref{LD1} and \eqref{LD2} follows from
    \eqref{BSMN} as   in the proof of Theorem \ref{T3.1}, together with the corresponding estimate on $\E|D^*(h^\aa)|^2$ as in the proof of \cite[Theorem 1.1]{WZ13}.  For instance, below we  outline the proof of \eqref{BSMN}.

 Firstly, for $s\in (0,1)$ let $X_t^s$ solve \eqref{E5} with $X_0^{\phi,s}= X_0+s\phi(X_0)$, let  $\mu^{\phi,s}= \L_{X_0^{\phi,s}}=\mu\circ({\rm Id}+\phi)^{-1},$
 and let $\aa_{r,t}^{\phi,s}$ be defined as $\aa_{r,t}$ with $X_t^{\phi,s}$ replacing $X_t$. Then as in \eqref{*DR} and \eqref{W001}, \eqref{INTN2} implies
\beq\label{TTS0} \beg{split} (P_Tf)(\mu^{\phi,\vv})- (P_T f)(\mu) &= \int_0^\vv   \E\<(\nn f)(X_T^{\phi,s}), \nn_{\phi(X_0)}X_T^{\phi,s}\>\,\d s\\
&= \int_0^\vv \E\big[f(X_T^{\phi,s}) D^*(h^{\aa^{\phi,s}})\big],  \ \ f\in C_b^1(\R^{m+d}),\end{split}\end{equation}
where $h^{\aa^{\phi,s}} := h_{0,\cdot}^{\aa^{\phi,s}_{0,\cdot}}$ satisfies
\beq\label{TTS1} \lim_{s\to 0} \E|D^*(h^{\aa^{\phi,s}})-D^*(h)|^2 =0.\end{equation}
By the argument leading to \eqref{WPP02}, \eqref{TTS0} yields
$$\ff{(P_Tf)(\mu^{\phi,\vv})-(P_Tf)(\mu)}\vv = \ff 1 \vv\int_0^\vv \E\big[f(X_T^{\phi,s}) D^*(h^{\aa^{\phi,s}})\big]\,\d s,\ \ f\in \B_b(\R^{m+d}).$$
  Combining this with \eqref{TTS1},   we prove \eqref{BSMN}  provided
  \beq\label{TTS2} \lim_{\vv\downarrow 0} \ff 1 \vv\int_0^\vv \E\big[\{f(X_T^{\phi,s}) -f(X_T)\} D^*(h^{\aa})\big]\,\d s =0.\end{equation}
     For any $r\in (0,T)$, let $I_r= \E(D^*(h^\aa)|\F_r)$. By \eqref{TTS0} we obtain
     \beg{align*} &  \E \big[ \{f(X_T^{\phi,\vv})-f(X_T)\} I_r\big] =    \E\big[I_r\E(f(X_T^{\phi,\vv})-f(X_T)|\F_r)\big]\\
     &=     \E\bigg[I_r \int_0^\vv \E\big(\<\nn f(X_T^{\phi,s}), \nn X_T^{\phi,s}\> \big|\F_r\big)\,\d s \bigg]
      =    \E\bigg[I_r \int_0^\vv \E\big(f(X_T^{\phi,s}) D^*(h_{r,\cdot}^{\aa_{r,\cdot}})  \big|\F_r\big)\,\d s \bigg]  \\
    &=    \int_0^\vv \E\big[I_r  f(X_T^{\phi,s}) D^*(h_{r,\cdot}^{\aa_{r,\cdot}})  \big]\,\d s,\ \ f\in C_b^1(\R^d).\end{align*}  Combining this with  the argument extending   \eqref{WPP02} from $f\in C_b^1(\R^d)$ to $f\in\B_b(\R^d)$, we obtain
     $$\E \big[ \{f(X_T^{\phi,\vv})-f(X_T)\} I_r\big] =      \int_0^\vv \E\big[I_r  f(X_T^{\phi,s}) D^*(h_{r,\cdot}^{\aa_{r,\cdot}})  \big]\,\d s,\ \ f\in \B_b(\R^d).$$ Consequently,
     $$\lim_{\vv\to 0} \E \big[ \{f(X_T^{\phi,\vv})-f(X_T)\} I_r\big] =0,\ \ f\in\B_b(\R^d), r\in (0,T).$$
    Then for any $r\in (0,T)$,
    \beg{align*} & \limsup_{\vv\downarrow 0} \bigg|\ff 1 \vv\int_0^\vv \E\big[\{f(X_T^{\phi,s}) -f(X_T)\} D^*(h^{\aa})\big]\,\d s \bigg| \\
    &=  \limsup_{\vv\downarrow 0} \bigg|\ff 1 \vv\int_0^\vv \E\big[\{f(X_T^{\phi,s}) -f(X_T)\} \cdot\{D^*(h^{\aa})-I_r\}\big]\,\d s \bigg|\\
    &\le 2\|f\|_\infty \E|D^*(h^\aa)- \E(D^*(h^\aa)|\F_r)|.\end{align*}
    Letting $r\uparrow T$ we derive \eqref{TTS2}, and hence prove \eqref{BSMN} as explained above.
    \end{proof}
\paragraph{Acknowledgement.} Financial support by the DFG through the CRC 1283 $``$Taming uncertainty and profiting from randomness and low regularity in analysis, stochastics and their applications'' is acknowledged. The authors would like to thank the referee and Professor Yulin Song for helpful comments.

\end{document}